\documentclass{amsart}
\usepackage{amssymb}
\usepackage{amsmath}
\usepackage{amsrefs}
\usepackage[colorlinks=true]{hyperref}
\usepackage{cleveref}
\usepackage{amscd}
\usepackage{amsthm}
\usepackage{mathtools}
\usepackage{tikz-cd}
\usepackage{mathrsfs}
\usepackage[margin=2.93cm]{geometry}

\allowdisplaybreaks

\newtheorem{theoremletter}{Theorem}

\newtheorem{theorem}{Theorem}[section]
\newtheorem{lemma}[theorem]{Lemma}
\newtheorem{proposition}[theorem]{Proposition}
\newtheorem{corollary}[theorem]{Corollary}

\newtheorem{conjecture}[theorem]{Conjecture}
\newtheorem{question}[theorem]{Question}

\theoremstyle{definition}
\newtheorem{definition}[theorem]{Definition}
\newtheorem{example}[theorem]{Example}

\theoremstyle{remark}
\newtheorem{remark}[theorem]{Remark}

\numberwithin{equation}{section}

\newcommand{\R}{\mathbb{R}}

\newcommand{\grad}{\nabla}
\newcommand{\N}{\mathbb{N}}
\newcommand{\Z}{\mathbb{Z}}

\renewcommand{\L}{\mathcal{L}}

\newcommand{\lt}{\left}
\newcommand{\rt}{\right}
\newcommand{\tms}{\times}

\newcommand{\rmk}{\begin{remark}}
\newcommand{\ermk}{\end{remark}}
\newcommand{\cor}{\begin{corollary}}
\newcommand{\ecor}{\end{corollary}}
\newcommand{\eq}{\begin{equation}}
\newcommand{\eeq}{\end{equation}}
\newcommand{\eqs}{\begin{equation*}}
\newcommand{\eeqs}{\end{equation*}}
\newcommand{\prop}{\begin{proposition}}
\newcommand{\eprop}{\end{proposition}}
\newcommand{\thm}{\begin{theorem}}
\newcommand{\ethm}{\end{theorem}}
\newcommand{\conj}{\begin{conjecture}}
\newcommand{\econj}{\end{conjecture}}
\newcommand{\lem}{\begin{lemma}}
\newcommand{\elem}{\end{lemma}}
\newcommand{\defi}{\begin{definition}}
\newcommand{\edefi}{\end{definition}}
\newcommand{\ex}{\begin{example}}
\newcommand{\eex}{\end{example}}
\newcommand{\alis}{\begin{align*}}
\newcommand{\ealis}{\end{align*}}
\newcommand{\pf}{\begin{proof}}
\newcommand{\epf}{\end{proof}}
\newcommand{\ali}{\begin{align}}
\newcommand{\eali}{\end{align}}
\newcommand{\qus}{\begin{question}}
\newcommand{\equs}{\end{question}}

\newcommand{\mc}{\mathcal}
\renewcommand{\bf}{\textbf}

\newcommand{\C}{\mathbb{C}}

\newcommand{\sub}{\subset}

\newcommand{\ov}{\overline}

\newcommand{\bb}{\mathbb}

\newcommand{\op}{\operatorname}

\renewcommand{\a}{\alpha}
\renewcommand{\b}{\beta}
\renewcommand{\d}{\partial}
\newcommand{\e}{\epsilon}

\newcommand{\g}{\gamma}

\newcommand{\s}{\sigma}
\renewcommand{\t}{\theta}

\renewcommand{\l}{\lambda}
\renewcommand{\o}{\omega}

\newcommand{\G}{\Gamma}
\renewcommand{\L}{\Lambda}

\renewcommand{\S}{\Sigma}

\newcommand{\Q}{\mathbb{Q}}

\renewcommand{\ov}{\overline}

\begin{document}
\title{On linking of Lagrangian tori in $\R^4$}
\author{Laurent C\^{o}t\'{e}}
\date{\today}
\begin{abstract} We prove some results about linking of Lagrangian tori in the symplectic vector space $(\R^4, \o)$. We show that certain enumerative counts of holomophic disks give useful information about linking. This enables us to prove, for example, that any two Clifford tori are unlinked in a strong sense. We extend work of Dimitroglou Rizell and Evans on linking of monotone Lagrangian tori to a class of non-monotone tori in $\R^4$ and also strengthen their conclusions in the monotone case in $\R^4$. \end{abstract}
\maketitle

\section{Introduction}  Let $L_1$ and $L_2$ be disjoint Lagrangian tori in the symplectic vector space $(\R^4, \o)$ where $\o= dx_1\wedge dy_1  + dx_2\wedge dy_2$. We say that $L_1$ and $L_2$ are \emph{smoothly unlinked} if they can be isotoped away from each other without intersecting. A more precise definition is as follows.

\defi \label{definition: smoothunlinking} Two closed, disjoint submanifolds $N_1, N_2 \sub \R^m$ for $m>1$ are said to be \emph{smoothly unlinked} if there exists a smooth isotopy $\phi^{(1)}: N_1 \tms [0,1] \to \R^m$ with $\phi_0^{(1)}= \op{Id}$ such that
\begin{itemize}
\item[(i)] $\phi_t^{(1)}(N_1) \cap N_2 = \emptyset$ for all $t \in [0,1]$,
\item[(ii)] $\phi_1^{(1)}(N_1)$ and $N_2$ are contained in disjoint, embedded balls. 
\end{itemize}

We say that $N_1$ and $N_2$ are \emph{smoothly linked} if they are not smoothly unlinked. By the isotopy extension theorem, the existence of $\phi^{(1)}$ is equivalent to the existence of an isotopy $\phi^{(2)}$ satisfying properties (i) and (ii) with the roles of $N_1$ and $N_2$ interchanged. \edefi

The following theorem is one of the main results of this paper. 

\begin{theoremletter} \label{theorem:clifforddisjoint} Let $L_1, L_2 \sub \R^4$ be disjoint Clifford tori of possibly different monotonicity factor. Then $L_1$ and $L_2$ are smoothly unlinked. \end{theoremletter}

For $r>0$, we say that $L \sub \R^4$ is a Clifford torus of monotonicity factor $\pi r^2/2$ if it is Hamiltonian isotopic to the standard model $\{(z_1, z_2) \in \C^2 \mid |z_1|=|z_2| = r\}$.  The proof of \Cref{theorem:clifforddisjoint} will be provided at the end of \Cref{section: enumerativeinvariant}; see \Cref{corollary:cliffordunlink1}.  

We remark that, in contrast to \Cref{theorem:clifforddisjoint}, any Lagrangian torus in $(\R^4, \o)$ is smoothly linked with a Chekanov torus; see \Cref{example:chekanovlink}.

Important progress in understanding linking of Lagrangian tori was achieved by Dimitroglou Rizell and Evans \cite{evans-rizell} using the theory of punctured pseudoholomorphic curves. Earlier work using different tools includes \cites{polterovich, eli-pol, borrelli}. 

In this paper, we build on ideas introduced in \cite{evans-rizell} and \cite{igr} to prove new results on linking in the restricted context of Lagrangian tori in $\R^4$. In particular, we prove the following theorem. 

\begin{theoremletter} \label{theorem:secondtheorem} Let $L_1, L_2 \sub (\R^4, \o)$ be disjoint, monotone Lagrangian tori with monotonicity factor $K_1$ and $K_2$ respectively. If $K_2 > K_1,$ then $L_1$ and $L_2$ are smoothly unlinked if and only if the image of the natural map $\pi_1(L_1) \to \pi_1(\R^4-L_2)$ vanishes. If $K_1=K_2$, then $L_1$ and $L_2$ are always smoothly unlinked. \end{theoremletter}

The proof of \Cref{theorem:secondtheorem} occupies most of \Cref{section:monotonelink}, where it is stated in a more general form as \Cref{corollary:restateB}

We emphasize that our proof of \Cref{theorem:secondtheorem} relies crucially on the special properties of holomorphic curves in dimension $4$. In contrast, the results of \cite{evans-rizell} work equally well in all dimensions; cf. \Cref{theorem:edr}.

The results of \cite{evans-rizell} on linking are restricted to monotone Lagrangian submanifolds. Our next theorem extends some of the results of \cite{evans-rizell} to a class of non-monotone Lagrangian tori in $\R^4$. In order to state precisely what is involved, we make a short digression to collect some necessary definitions.

Given a symplectic manifold $(M, \o)$ and a Lagrangian submanifold $L \sub M$, the \emph{Maslov class} is a map $ \mu: \pi_2(M, L) \to \Z $ which takes values in the even integers if $L$ is orientable. By a slight abuse of notation, the \emph{symplectic area class} is defined as the map $\o: \pi_2(M, L) \to \R$ taking $[u] \mapsto \o([u])= \int_{D^2} u^* \o.$

The following invariant of Lagrangian tori in $\R^4$ will play an essential role throughout this paper. 

\defi \label{definition:a2} Let $L \sub (\R^4, \o)$ be a Lagrangian torus. We define 
\eqs A_2(L):= \op{min} \{ \o(\a) \mid \a \in \pi_2(\R^4, L),\, \mu(\a)=2, \, \o(\a)>0 \}. \eeqs
\edefi

The following definition was considered in \cite{evans-rizell}, using slightly different terminology.

\defi[cf. \cite{evans-rizell}] \label{definition: ERlinking} Let $N_1$ and $N_2$ be closed, disjoint submanifolds of $\R^{m}$.  Then $N_1$ is said to be \textit{homologically unlinked} from $N_2$ if $[N_1] \in H_2(\R^{m}-N_2; \Z)$ is the zero class. Otherwise, we say that $N_1$ is \emph{homologically linked} with $N_2.$ We say that $N_1$ and $N_2$ are \textit{mutually homologically unlinked} if each one is null-homologous in the complement of the other. \edefi

Clearly smooth unlinking implies homological unlinking. Observe also that the notion of homological linking is not symmetric, i.e. it may be the case  that $N_1$ is homologically unlinked from $N_2$ while $N_2$ is homologically linked with $N_1$; see \Cref{example:notsmoothlydisjoin}. This is in contrast to the notion of smooth unlinking (see \Cref{definition: smoothunlinking}), which is manifestly symmetric. 

In \Cref{section: non-monotone}, we introduce a class of non-monotone Lagrangian tori called \emph{admissible}. These tori are distinguished by the nonvanishing of an enumerative invariant which counts Maslov $2$ disks of small area; cf. \Cref{definition:admissible}. We show that the class of admissible tori is closed under Hamiltonian isotopy, and that it contains ``most" product tori. 

As mentioned above, the results of \cite{evans-rizell} on linking only concern monotone Lagrangian submanifolds. The following theorem extends \cite[Theorem A]{evans-rizell} to admissible tori in $\R^4$. 

\begin{theoremletter} \label{theorem: nonmonotone} Let $L_1, L_2 \sub \R^4$ be disjoint Lagrangian tori and suppose that $L_1$ is admissible. If $A_2(L_2) \geq A_2(L_1),$ then $[L_1]$ is the zero class in $H_2(\R^4-L_2; \Z).$ In other words, $L_1$ is homologically unlinked from $L_2$. \end{theoremletter}

In \Cref{section: linkingconstruction}, we show that the assumption $A_2(L_2) \geq A_2(L_1)$ in \Cref{theorem: nonmonotone} is sharp in a suitable sense; see \Cref{proposition: necessary1}.

\subsection{Some perspective} \label{subsection:comparison}  One of the main conclusions of this paper may be summarized as follows: if one considers the problem of smoothly unlinking monotone Lagrangian tori in $\R^4$, then the obvious algebro-topological obstructions are the only obstructions. Moreover, one can identify reasonable conditions under which these obstructions vanish.

Suppose that $L_1, L_2 \sub (\R^4, \o)$ are monotone Lagrangian tori. In order for $L_1$ and $L_2$ to be smoothly unlinked, it is necessary that the natural maps $H_k(L_i; \Z) \to H_k(\R^4-L_j; \Z)$ and $\pi_k(L_i) \to \pi_k(\R^4-L_j)$ have trivial image for all $k \geq 0$ and $i \neq j \in \{1,2\}$. In addition to these algebro-topological obstructions, there could a priori be more subtle obstructions coming from smooth topology. Indeed, there is in general a large gap between algebraic and differential topology in dimension 4. 

Let us now assume without loss of generality that the monotonicity factor of $L_1$ is at most equal to that of $L_2$. Under this assumption, we will show in \Cref{section:monotonelink} that $L_1$ bounds a solid torus which is smoothly embedded in the complement of $L_2$; see \Cref{theorem: erimprove}. In particular, this implies that $L_1$ is homologically unlinked from $L_2$, which was already proved by Dimitroglou Rizell and Evans; see \Cref{theorem:edr}. 

As  noted above, it is necessary, in order for $L_1$ and $L_2$ to be smoothly unlinked, that the natural map $\pi_1(L_1) \to \pi_1(\R^4-L_2)$ have trivial image. Using the fact that $L_1$ bounds a smoothly embedded solid torus, we show that this necessary condition is in fact \emph{sufficient} (see \Cref{theorem:secondtheorem} and \Cref{corollary:restateB}).  Hence the question of whether $L_1$ and $L_2$ are smoothly unlinked reduces to elementary algebraic topology.

In \Cref{section: enumerativeinvariant}, we analyze the map $\pi_1(L_1) \to \pi_1(\R^4-L_2)$.  We show that it must have trivial image if certain enumerative counts of holomorphic disks with boundary in $L_1$ are nonzero.  This enables us, in particular, to prove that Clifford tori (of possibly different monotonicity factors) are always smoothly unlinked; see \Cref{theorem:clifforddisjoint}.

\subsection{Organization}  \Cref{section:prep} contains a summary of some prior work on linking of Lagrangian tori, and some topological lemmas which will be needed in the remainder of the paper. 

\Cref{section:monotonelink} and \Cref{section: enumerativeinvariant} were already surveyed in the above paragraphs; they contain in particular proofs of \Cref{theorem:clifforddisjoint} and \Cref{theorem:secondtheorem}. 

\Cref{section: non-monotone} deals with linking of non-monotone tori in $\R^4$. In particular, we introduce the class of admissible tori alluded to earlier and prove \Cref{theorem: nonmonotone}. 

\Cref{section: linkingconstruction} describes a construction which shows that \Cref{theorem: nonmonotone} is sharp in a suitable sense. 

\Cref{section:quantitativeunlink} explores some connections between our analysis of linking and questions about embeddings of tori and polydisks into various subdomains of $\R^4$.

\bf{Acknowledgements. } I wish to particularly thank my advisor Yasha Eliashberg for his guidance throughout this project and for many decisive insights and suggestions.  I also wish to thank Kai Cieliebak for suggestions which significantly strengthened \Cref{section:quantitativeunlink} of this paper. I wish to thank Georgios Dimitroglou Rizell for kindly explaining several aspects of \cite{igr} to me, for suggesting the argument in \Cref{remark:charette} and for several other useful comments on a draft of this paper. During the course of this project, I benefited from conversations with Daniel \'{A}lvarez-Gavela, Tobias Ekholm, C\'{e}dric de Groote, Eleny Ionel, Janko Latschev, Oleg Lazarev, Daniel Ruberman, Laura Starkston and Chris Wendl. I gratefully thank all of them. Finally, I wish to thank the anonymous referee for many important comments and suggestions.

\section{Context and preparatory material}  \label{section:prep} This section is intended to introduce some preparatory material and to provide some context to help motivate the techniques and results of this paper. We begin by stating some standard conventions which will be followed throughout this work. We then summarize some prior work on linking of Lagrangian submanifolds. Finally, we prove some topological lemmas which will be useful in later sections and which partly rely on an important classification theorem of Dimitroglou Rizell, Ivrii and Goodman.

\subsection{Conventions} \label{subsection:conventions}
Unless otherwise indicated, the vector space $\R^{2n}$ is endowed with the coordinates $(x_1,y_1,\dots,x_n,y_n)$ and with the symplectic form $\o= dx_1 \wedge dy_1 +\dots+ dx_n \wedge dy_n.$ We let $j$ denote the standard integrable complex structure on $\R^{2n}$. We will routinely identify $\R^{2n}$ with $\C^n$ via the map $(x_1,y_1,\dots,x_n,y_n) \mapsto (x_1+ iy_1,\dots, x_n+ i y_n).$ 

Given a Lagrangian $L \sub (\R^{2n}, \o)$, we note that the boundary maps $\pi_*(\R^{2n}, L) \to \pi_{*-1}(L)$ and $H_*(\R^{2n}, L; \Z) \to H_{*-1}(\R^4, L; \Z)$ are isomorphisms. It follows from the universal coefficient theorem that we may view the Maslov class $\mu$ and symplectic area class $\o$ as cohomology classes of $L.$

\defi \label{definition:torus} Given a symplectic manifold $(M, \o)$, a Lagrangian submanifold $L \sub (M, \o)$ is said to be a \emph{Lagrangian torus} if it is diffeomorphic to $\bb{T}^n= S^1\tms\dots\tms S^1$. \edefi

\defi \label{definition: monotone} Given a symplectic manifold $(M, \o)$, a Lagrangian submanifold $L \sub (M, \o) $ is said to be \emph{monotone} if $\o(\a)= c \mu(\a)$ for all $\a \in \pi_2(M, L)$. Here $c$ is a positive constant which is called the \emph{monotonicity factor} of $L$. \edefi

For a Lagrangian torus $L \sub (\R^4, \o)$, note that the monotonicity factor $c$ satisfies the identity ${A_2(L)}= 2c$; cf. \Cref{definition:a2}. 

\ex \label{example: clifford} For $r>0$, the \emph{Clifford torus} \eqs L_{Cl}= \{ (z_1, z_2) \in \C^2 \mid |z_1|= |z_2|=r \} \eeqs is a monotone Lagrangian torus of monotonicity factor $\pi r^2/2$.  \eex

\ex \label{example: chekanov} For $r>0$, the \emph{Chekanov torus} is defined as the set \eqs L_{Ch} = \{  ((e^x+i e^{-x}y) \cos \t, (e^x+i e^{-x}y) \sin \t) \mid \t \in [0, 2\pi],  x^2+y^2= r^2, (x, y) \in \R^2\}. \eeqs
It is a monotone Lagrangian torus of monotonicity factor $\pi r^2/2$. \eex

A Lagrangian torus in $\R^{4}$ which is Hamiltonian isotopic to $L_{Cl}$ or $L_{Ch}$ for some $r>0$ will be referred to as a Clifford torus or a Chekanov torus. When it is not clear from the context, we will indicate whether we are referring the standard models of \Cref{example: clifford} and \Cref{example: chekanov} or to a torus Hamiltonian isotopic to them.

A \emph{pseudoholomorphic} or \emph{J-holomorphic} curve is a map $u: (\S, j) \to (W^{2n}, J)$ satisfying the nonlinear Cauchy-Riemann equation $du+ J \circ du \circ j=0$. Here $(\S, j)$ is a (possibly punctured) Riemann surface and $(W^{2n}, J)$ is an almost-complex manifold. Such maps will routinely be referred to as \emph{holomorphic curves} when it is clear from the context that $J$ is not assumed to be integrable. 

By a similar abuse of language, we will usually use the terms \emph{almost-complex structure} and \emph{complex structure} interchangeably, even though the later term is often reserved in the literature for integrable almost-complex structures.   

\rmk[Signs] \label{remark:signs} We will generally follow the sign conventions of \cite{evans-rizell} and \cite{igr}.  In particular, the Liouville $1$-form on a cotangent bundle $T^*M$ is denoted $\l$ and gives rise to a symplectic form $\o$ by the equation $\o= d\l$. We note that this sign convention differs from that of \cite[see Remark 3.5.35]{mcduff-sal-intro}. \ermk

\subsection{Some prior work}  \label{subsection:previouswork} Let $\phi: M \to (\R^{2n}, j)$ be a totally real embedding. Given a nowhere vanishing vector field $X \in \G(TM),$ let $M'$ be a small push-off of $\phi(M)$ in the direction of $j(d\phi(X)).$ This gives rise to a class $[M'] \in H_n(\R^{2n}- \phi(M); \Z).$ By the long exact sequence of the pair $(\R^{2n}, \R^{2n}-\phi(M))$ and Alexander duality, there are isomorphisms 
\eq \label{equation: theisom} H_n(\R^{2n} - \phi(M); \Z) \simeq H_{n+1}(\R^{2n}, \R^{2n}- \phi(M); \Z) \simeq H^{n-1}(M; \Z). \eeq

\defi Let $l(\phi, X) \in H^{n-1}(M; \Z)$ be the class corresponding to $[M']$ under \eqref{equation: theisom}. The class $l(\phi, X)$ is called the \emph{linking class}. \edefi

Observe that the linking class $l(\phi, X)$ vanishes if and only if $M'$ is homologically unlinked from $\phi(M)$; cf. \Cref{definition: ERlinking}.

It can be shown that there exist totally real embeddings $\phi: \bb{T}^2 \to \R^4$ and vector fields $X \in \G(T\bb{T}^2)$ such that $l(\phi, X) \neq 0.$ In contrast, for Lagrangian embeddings Eliashberg and Polterovich proved the following theorem using the technique of Luttinger surgery.

\thm[Eliashberg--Polterovich \cite{eli-pol}] \label{theorem: eliash-polt} Let $i: \bb{T}^2 \to (\R^4, \o)$ be a Lagrangian embedding. Then $l(i, X)=0$ for all nonvanishing vector fields $X \in \G(TM).$ \ethm
We remark that \Cref{theorem: eliash-polt} was extended by Borrelli \cite{borrelli} to Lagrangian embeddings of $S^1 \tms S^3$ and $S^1 \tms S^7$ into $\R^8$ and $\R^{16}$ respectively. 

In \cite{evans-rizell}, Dimitroglou Rizell and Evans introduced a new approach to the study of linking of Lagrangian submanifolds. This approach relies on the theory of punctured pseudoholomophic curves. Dimitroglou Rizell and Evans proved the following theorem.  

\thm[Dimitroglou Rizell--Evans, Theorem A in \cite{evans-rizell}] \label{theorem:edr} Let $W$ be a subcritical Stein manifold and let $K_2 \geq K_1>0$ be real numbers. An embedded monotone Lagrangian torus $L_1$ with monotonicity factor $K_1$ is homologically unlinked from any embedded monotone Lagrangian torus $L_2$ with factor $K_2$. In particular, two embedded monotone Lagrangian tori with the same monotonicity factor are mutually homologically unlinked. \ethm

Let us briefly sketch the proof of \Cref{theorem:edr} as it is the starting point for much of this work. We refer the reader to \cite{evans-rizell} for details. 

Given a class $\b \in \pi_2(W, L_1)$ and an almost-complex structure $J,$ let $\mc{M}_{1}(\b, J)$ be the moduli space of $J$-holomorphic disks representing the class $\b$ with one interior marked point. Let $\mc{M}_{0,1}(\b, J)$ be the corresponding moduli space of $J$-holomorphic disks with one boundary marked point.  Under suitable assumptions on $J,$ it can be shown using work of Damian \cite{damian} and Evans-K\c{e}dra \cite{ev-k} that there exists a class $\b$ such that the boundary evaluation map $\mc{M}_{0,1}(\b, J) \to L_1$ has nonzero degree on some component $\mc{M}$ of the moduli space $\mc{M}_{0,1}(\b, J)$. This fact relies on the assumption that $L_1$ is monotone.

The crux of the argument is now to produce an almost-complex structure $J$ with the property that the image of $\mc{M}_{1}(\b, J)$ under the natural evaluation map is disjoint from $L_2.$  This can be achieved by deforming a fixed complex structure $J_0$ near $L_2$ by a process known as ``stretching the neck". The authors analyze the limiting behavior of sequences of holomorphic disks under the SFT compactness theorem. Using the assumption that $K_2 \geq K_1,$ they conclude that all disks must become disjoint from $L_2$ for sufficiently large deformations of the complex structure. The theorem then follows by elementary topological arguments since $\d [\mc{M}]= n[L_1]$ for some $n \geq 1$.

We remark that arguments similar to the one sketched above appear in a recent paper of Ekholm and Smith \cite[see Thm. 1.3]{ekholm-smith}.

\subsection{Topological lemmas} We now state a landmark classification result of Ivrii, Goodman and Dimitroglou Rizell. Both the theorem and its proof will have an important role in our work.

\thm[Dimitroglou Rizell--Ivrii--Goodman, \cite{igr}] \label{theorem: igr} All Lagrangian tori in $(\R^4, \o),$ $(S^2 \tms S^2, \o \oplus \o)$ and $(\C P^2, \o_{\op{FS}})$ are isotopic through Lagrangian tori. \ethm

One can show by elementary topological arguments that all orientable Lagrangian submanifolds of $\R^4$ and $\C P^2$ are tori, and that all orientable Lagrangian submanifolds of $S^2 \tms S^2$ are spheres or tori. Hence \Cref{theorem: igr} gives a complete classification of all Lagrangian submanifolds of $\R^4$ and $\C P^2$ up to Lagrangian isotopy. One also obtains a full classification of Lagrangian submanifolds of $S^2 \tms S^2$, up to Lagrangian isotopy, by combining \Cref{theorem: igr} with a theorem of Hind \cite{hind} establishing the uniqueness of Lagrangian spheres in $S^2 \tms S^2$ up to Hamiltonian isotopy. 

We now introduce some useful topological lemmas whose proofs rely on \Cref{theorem: igr}. 

\lem \label{lemma: homotopyofcomplement} Let $L \sub \R^4$ be a Lagrangian torus.   Then $\pi_1(\R^4-L) = H_1(\R^4-L; \Z) = \Z.$  Moreover, we have that $H_2(\R^4-L; \Z)= \Z \oplus \Z$ and $H_3(\R^4-L; \Z)= \Z.$  For all $i \geq 4,$ the groups $H_i(\R^4-L; \Z)$ vanish. 
\elem

\pf It follows from \Cref{theorem: igr} that all Lagrangian tori in $\R^4$ are Lagrangian isotopic. Hence we may assume that $L$ is the Clifford torus of radius one $\{ (z_1, z_2) \in \C^2 \mid |z_1|= |z_2|=1 \}.$  Let $\mc{U} := \{(x_1, y_1, x_2, y_2) \in \R^4 \mid x_1^2 +y_1^2 \neq 0 \}.$ One can easily check (e.g. using cylindrical coordinates) that $\pi_1(\mc{U}-L)=\Z  \oplus \Z.$  Let $\mc{V}:= \{ (x_1, y_1, x_2, y_2) \mid x_1^2+ y_1^2 < 1/2\}.$  An application of van-Kampen's theorem implies that $\pi_1(\R^4-L)=  \pi_1( \mc{U}- L) *_{\pi_1(\mc{U} \cap \mc{V})} \pi_1(\mc{V}) =(\Z \oplus \Z) *_{\Z} \{ e\} = \Z.$   

To compute the homology groups, let $N(L)$ be a tubular neighborhood of $L$ and consider the Mayer-Vietoris homology sequence associated to the subspaces $N(L) \sub \R^4$ and $\R^4-L \sub \R^4.$  Observing that $N(L) \cap (\R^4-L)$ is homotopy equivalent to $\bb{T}^3,$ we find that $H_k(\bb{T}^3) = H_k(L) \oplus H_k(\R^4-L).$ Since $H_2(\bb{T}^3)= \Z \oplus \Z \oplus \Z$ and $H_3(\bb{T}^3)= \Z,$ it follows easily that $H_2(\R^4-L)= \Z \oplus \Z$ and $H_3(\R^4-L)= \Z.$  The vanishing of the higher homology groups holds for dimension reasons. \epf

\lem \label{lemma: pi1link} Let $L_1$ and $L_2$ be disjoint Lagrangian tori in $\R^4.$ If the natural map $i_*: \pi_1(L_1) \to \pi_1(\R^4-L_2)$ has nontrivial image, then $L_2$ is homologically linked with $L_1$; cf. \Cref{definition: ERlinking}. \elem

\pf We first argue that there exists a solid torus $S \sub \R^4$ with the property that $\d S=L_2$ and that the intersection pairing
\begin{align} \label{equation:intersectionpairing} 
H_3(\R^4, L_2; \Z) \tms H_1(\R^4-L_2; \Z) &\to \Z \\
([S], [\g]) &\mapsto S \cdot \g \nonumber 
\end{align}
generates $\op{Hom}(H_1(\R^4-L_2; \Z), \Z) \simeq H^1(\R^4-L_2; \Z).$

This is not hard to verify in the case where $L_2$ is a Clifford torus (in this case, there are two families of disks which form solid tori with the desired property).  \Cref{theorem: igr} implies that $L_2$ is isotopic to the Clifford torus, so the claim follows in general by the isotopy extension theorem. 

Let us now assume for contradiction that $i_*: \pi_1(L_1) \to \pi_1(\R^4-L_2)$ has nontrivial image and that $[L_2] \in H_2(\R^4-L_1; \Z)$ is the zero class. Let $[\g] \in \op{Im} i_*$ be a nonzero element. In light of \Cref{lemma: homotopyofcomplement}, we have isomorphisms $\Z \simeq \pi_1(\R^4-L_2) \simeq H_1(\R^4-L_2; \Z)$, so we may as well view $[\g]$ as a nonzero class in $H_1(\R^4-L_2; \Z).$

Since $0=[L_2] \in H_2(\R^4-L_1; \Z)$, there exists some chain $U \in C_3(\R^4-L_1)$ such that $L_2= \d U$. Observe that $U$ naturally defines a class in $H_3(\R^4, L_2; \Z)$, and we have $U \cdot [\g]=0$ since $U$ is disjoint from $L_1$. 

It follows from the long exact sequence of the pair $(\R^4, L_2)$ that $0=(U-[S]) \in H_3(\R^4, L_2; \Z)=\Z. $  Hence $0=(U-[S]) \cdot [\g] = U \cdot [\g] - [S] \cdot [\g]$.  Hence $[S] \cdot [\g]=0$. In light of \eqref{equation:intersectionpairing}, this implies that $[\g]=0$ as an element of $H_1(\R^4-L_2; \Z).$ We thus obtain a contradiction. \epf

The converse of \Cref{lemma: pi1link} is true under the assumption that $L_1$ and $L_2$ are monotone Lagrangian tori in $\R^4$. This will follow from \Cref{corollary: pi1linking}. 

We end this section with an example which was already mentioned in the introduction. Given disjoint compact Lagrangians $L_1, L_2 \sub \R^{2n}$, this example illustrates that $L_1$ may be homologically linked with $L_2$ while $L_2$ is homologically unlinked from $L_1$. 

\ex \label{example:notsmoothlydisjoin} Let $L_1 = \{ (z_1, z_2) \in \C^2 \mid |z_1|=|z_2|=1 \}$ be the Clifford torus.   We showed in \Cref{lemma: homotopyofcomplement} that $\pi_1(\R^4-L_1)= \Z.$  Choose a loop $\g$ realizing a nontrivial element of $\pi_1(\R^4-L_1).$   It follows easily from the isotropic neighborhood theorem that any arbitrarily small neighborhood of $\g$ contains a Lagrangian torus which is a circle bundle over $\g.$  

Let $L_2$ be such a Lagrangian torus. It follows by construction that $L_2$ is null-homologous in the complement of $L_1.$  However, the inclusion $L_2 \hookrightarrow \R^4-L_1$ induces a nontrivial map on fundamental groups. It now follows by \Cref{lemma: pi1link} that $L_1$ is homologically essential in the complement of $L_2.$  \eex

\section{Linking of monotone Lagrangian tori} \label{section:monotonelink}

The main result of this section is \Cref{theorem: erimprove}, which leads almost immediately to a proof of \Cref{theorem:secondtheorem} (cf. \Cref{corollary:restateB}) and is also an essential ingredient in the proof of \Cref{theorem:clifforddisjoint}. The arguments of this section borrow heavily from \cite{igr} and \cite{evans-rizell}, but we have included most proofs since our setting is slightly different. 

\subsection{The main result} We recall from the introduction the following definition, which plays an essential role throughout this work. 

\defi Let $L \sub (\R^4, \o)$ be a Lagrangian torus. Then 
\eqs A_2(L):=  \op{min} \{ \o(\a) \mid \a \in \pi_2(\R^4, L),\, \mu(\a)=2, \, \o(\a)>0 \}. \eeqs
\edefi

The goal of this section is to prove the following theorem. 

\thm \label{theorem: erimprove} Let $L_1, L_2 \sub (\R^4, \o)$ be disjoint Lagrangian tori.  Assume that $L_1$ is monotone and that $A_2(L_2) \geq A_2(L_1).$  Then there exists a smooth embedding $$\phi: (S^1 \tms D^2, S^1 \tms \d D^2) \to (\R^4-L_2, L_1).$$  \ethm

In other words, \Cref{theorem: erimprove} says that $L_1$ bounds a solid torus in the complement of $L_2$.  

\cor \label{corollary: pi1linking} Suppose that $L_1^{(1)}, L_1^{(2)}, L_2$ are Lagrangian tori in $(\R^4, \o)$ with $L_1^{(i)}$ monotone and with $L_1^{(i)} \cap L_2 = \emptyset$ for $i=1,2.$  If $A_2(L_2) \geq A_2(L_1^{(i)})$, then $L_1^{(1)}$ and $L_1^{(2)}$ are smoothly isotopic in the complement of $L_2$ if and only if the group homomorphisms $\pi_1(L_1^{(i)}) \to \pi_1(\R^4-L_2)$ have the same image. \ecor

\pf[Proof of \Cref{corollary: pi1linking}] If $L_1^{(1)}$ and $L_1^{(2)}$ are smoothly isotopic in $\R^4-L_2,$ then clearly the images of the induced maps of fundamental groups coincide. For the reverse direction, observe by the theorem that $L_1^{(i)}$ bounds a solid torus $\phi^{(i)}: (S^1 \tms D^2, S^1 \tms \d D^2) \to (\R^4-L_2, L_1^{(i)}).$ 

The interior of $\op{Im} \phi^{(i)}$ forms an orientable open submanifold and hence has trivial normal bundle. It follows that there exists for some $\e>0$ a tubular neighborhood of the circle $\Psi^{(i)}: S^1 \tms \R^2 \tms \R \to \R^4$ such that $\Psi^{(i)}(t, x_1, x_2, 0)= \phi^{(i)}(t, x_1, x_2)$ for $|x|< 2\e$. 

Let $T^{(i)}= \{ \Psi^{(i)}(t, \e \cos \t, \e \sin \t, 0) \mid t \in S^1, \t \in [0, 2\pi)\}$.  By contracting $L_1^{(i)}$ using $\phi^{(i)}$, we immediately see that $L_1^{(i)}$ and $T^{(i)}$ are isotopic in the complement of $L_2$. 

By hypothesis, the cores $\Psi^{(i)}(S^1 \tms 0 \tms 0)$ are isotopic in $\R^4-L_2$ for $i=1,2$. Since any isotopy between these cores extends to an isotopy of tubular neighborhoods, it follows (after perhaps choosing $\e$ smaller) that $\Psi^{(2)}$ is isotopic to some tubular neigborhood $\tilde{\Psi}^{(2)}: S^1 \tms \R^2 \tms \R$, where $\Psi^{(1)}(S^1 \tms 0 \tms 0 \tms 0) = \tilde{\Psi}^{(2)}(S^1 \tms 0 \tms 0 \tms 0)$. Similarly, $T^2$ is isotopic to $\tilde{T}^2=  \{ \tilde{\Psi}^{(2)}(t, \e \cos \t, \e \sin \t, 0) \mid t \in S^1, \t \in [0, 2\pi)\}$.

The uniqueness theorem for tubular neighborhoods \cite[p. 112]{hirsch} implies that there exists a smooth isotopy of tubular neighborhoods $\Psi_{t}^{(1)}: [0,1] \tms S^1 \tms \R^2 \tms \R  \to \R^4-L_2$ with $\Psi_0^{(1)}= \Psi^{(1)}$ and such that $\Psi^{(1)}_1= \tilde{\Psi}^{(2)} \circ F$, where $F: S^1 \tms \R^3 \to S^1 \tms \R^3$ is a bundle isomorphism. 

After a possible further isotopy, we can assume that this bundle isomorphism is a fiberwise isometry, with respect to the standard euclidean metric on $\R^3$. Finally, we can assume that it preserves the splitting $\R^2 \tms \R$, since we can generate $\pi_1(SO(3))= \Z/2$ by a loop of orthogonal matrices which rotates the plane $\R^2 \tms 0$ around the axis $0 \tms 0 \tms \R$.  It follows that $T^{(1)}$ and $\tilde{T}^{(2)}$ are isotopic.  

Since $T^1$ is isotopic to $L_1^{(1)}$ and since $\tilde{T}^2$ is isotopic to $L_1^{(2)}$, it follows that $L_1^{(1)}$ are $L_1^{(2)}$ isotopic. \epf

The following corollary of \Cref{theorem: erimprove} was already stated in the introduction in a slightly weaker form as \Cref{theorem:secondtheorem}. It strengthens the conclusions of \Cref{theorem:edr}, due to Dimitroglou Rizell and Evans, in the special case of Lagrangian tori in $\R^4$.   

\cor \label{corollary:restateB} Let $L_1, L_2 \sub (\R^4, \o)$ be disjoint Lagrangian tori.  If $L_1$ is monotone and $A_2(L_2) \geq A_2(L_1)$, then $L_1$ bounds a solid torus in the complement of $L_2$. Moreover, $L_1$ and $L_2$ are smoothly unlinked if and only if the image of the natural map $\pi_1(L_1) \to \pi_1(\R^4-L_2)$ vanishes.  If $L_1$ and $L_2$ are both monotone and $A_2(L_1)= A_2(L_2)$, then $L_1$ and $L_2$ are smoothly unlinked. \ecor

\pf[Proof of \Cref{corollary:restateB}] The fact that $L_1$ bounds a solid torus in the complement of $L_2$ is a restatement of the theorem.

If $L_1$ and $L_2$ are smoothly unlinked, then it follows immediately from \Cref{definition: smoothunlinking} that the map $\pi_1(L_1) \to \pi_1(\R^4-L_2)$ has trivial image. To prove the converse, consider some other Lagrangian torus $L_1'$ which is far away from $L_2$ and, in particular, is smoothly unlinked from $L_2$. Then it follows from \Cref{corollary: pi1linking} that $L_1$ and $L_1'$ are isotopic in the complement of $L_2$ if the map $\pi_1(L_1) \to \pi_1(\R^4-L_2)$ has trivial image. Hence $L_1$ and $L_2$ are smoothly unlinked. 

In the special case where $A_2(L_1)=A_2(L_2)$, we find by interchanging the roles of $L_1$ and $L_2$ that they both bound solid tori in the complement of the other. It then follows from \Cref{lemma: pi1link} that the natural map $\pi_1(L_1) \to \pi_1(\R^4-L_2)$ has trivial image. Hence we conclude that $L_1$ and $L_2$ are smoothly unlinked.\epf

\pf[Overview of the proof of \Cref{theorem: erimprove}] \phantom\qedhere Our proof of \Cref{theorem: erimprove} is very much analogous to the original argument of Dimitroglou Rizell and Evans in their proof of \Cref{theorem:edr}. The main difference is that we work with holomorphic \emph{planes} rather than holomorphic disks. 

The argument of Dimitroglou Rizell and Evans was already sketched in the introduction; cf. \Cref{subsection:previouswork}. Given $L_1$ and $L_2$ as in the statement of \Cref{theorem:edr}, recall that the strategy is to deform the complex structure near $L_2$ by ``stretching the neck". One then considers the effect of this deformation on the moduli space of holomorphic disks with boundary on $L_1$. For the argument to work, one needs to ensure that the relevant moduli spaces remain non-empty as one deforms the complex structure. In the original paper of \cite{evans-rizell}, this property is observed to follow from work of Damian \cite{damian} and Evans-K\c{e}dra \cite{ev-k} using Floer theory. 

These Floer theoretic methods do not apply if one works with planes instead of disks. Instead, we will appeal to an analysis carried out in \cite{igr}, which also uses the technique of ``neck stretching" to produce moduli spaces of planes whose compactification has boundary in $L_1$. The relevant statement is \Cref{proposition: igrciting}. We will then analyze the behavior of these moduli spaces under deformation of the complex structure near $L_2$. This step is carried out in \Cref{proposition: sim1}. The argument and conclusion will be essentially the same as in \cite{evans-rizell}, although one can slightly sharpen the analysis when working in dimension $4$. This in particular allows us to replace the monotonicity assumption on $L_2$ with a condition on $A_2(L_2)$.

The monotonicity assumption on $L_1$ is needed in order to control the area of the holomorphic planes obtained using \cite{igr}. In \Cref{section: non-monotone}, we will prove certain results on homological linking of non-monotone Lagrangian tori in $\R^4$ using moduli spaces of holomorphic disks. It would be interesting to extend the arguments of that section to planes, but the analysis required seems more difficult; cf. \Cref{subsection:planesvsdisks}. \epf

\subsection{Recollection of some standard constructions} \label{subsection:neckstretch} For completeness and for the purpose of fixing some conventions which will be needed in the remainder of this work, we review some standard constructions on the way to proving \Cref{theorem: erimprove}. 
\defi Given a metric $g$ on $T^*M$ and a real number $r>0$, let $S^*_{r, g}M$ denote the sphere bundle consisting of covectors of norm $r$. Let $0_{M} \hookrightarrow T^*M$ denote the zero section. The submanifold $S^*_{r,g}M$ naturally inherits a contact structure $\a_{r,g}$ by restricting the Liouville form $\l_{can}$. \edefi

For $R > 0$, we consider the polydisk $\mc{P}(R, R)= \{ (z_1, z_2)  \in \C^2 \mid |z_1|< R, |z_2|<R \}$. We will be viewing $\mc{P}(R,R)$ both as an open symplectic manifold and as a symplectic subdomain of $(\R^4, \o)$. By choosing $R$ large enough, we can assume that $L_1$ and $L_2$ are both contained in $\mc{P}(R, R) \sub \mathbb{R}^4$. 

Observe that there is a natural symplectic embedding 
\eq i: \mc{P}(R, R) \to (S^2 \tms S^2, \o_R \oplus \o_R), \eeq 
where $\int \o_R= \pi R^2$ and $S^2 \tms S^2 = \mc{P}(R, R) \cup D_{\infty}$ with $D_{\infty}= S^2 \tms \{\infty\} \cup \{ \infty \} \tms S^2.$  Thus we may view $L_1$ and $L_2$ as Lagrangian submanifolds of $(S^2 \tms S^2, \o_R \oplus \o_R)$. 

It will be useful to consider the identification $T^*\bb{T}^2 \simeq \bb{T}^2 \tms \R^2$ given by the map $y_1 {d \t_1}+ y_2 d \t_2 \mapsto (\t_1, \t_2, y_1, y_2)$. In these coordinates, we have $\o_{can}= dy_1 \wedge d\t_1 + dy_2 \wedge d\t_2 $ and $\l_{can}= y_1 d\t_1+ y_2 d\t_2$; cf. \Cref{remark:signs}.

For $i=1,2$, let $\phi_i: \op{Op}(0_{\bb{T}^2}) \to N(L_i) \sub S^2 \tms S^2 - D_{\infty}$ be Weinstein embeddings with disjoint images. Let $g_i$ be a suitable rescaling of the flat metric on $T^*\bb{T}^2$ so that 
\eq \phi^i: \lt( (-1,1) \tms S^*_{1, g_i}\bb{T}^2, d(e^t\a_i) \rt) \to (N(L_i), \o_R \oplus \o_R) \eeq is a symplectic embedding. We write $\a_i= \a_{1, g_i}$. By setting $y_1= r \cos \ov{\t}$ and $y_2= r \sin \ov{\t}$ for $r>0$ and $\ov{\t} \in \R/\Z$, we naturally get coordinates $(\t_1, \t_2, \ov{\t}) \in (\R/\Z)^3$ on $S^*_{1,g_i} \bb{T}^2 \simeq \bb{T}^3$. We then have $\a_i= \e_i( \cos \ov{\t} d\t_1 + \sin \ov{\t} d\t_2)$ for some $\e_i>0$. Observe that the Reeb vector field is then $R_{\a_i}= \frac{1}{\e_i} (\cos \ov{\t} \d_{\t_1} + \sin \ov{\t} \d_{\t_2}).$

Let us consider the symplectization $(\R \tms S^*_{1,g_i}\bb{T}^2, d(e^t\a_i))$ with coordinates $(t, \t_1, \t_2, \ov{\t})$. We fix a trivialization of the tangent bundle $\Phi^i= \{ \d_t, R_{\a_i}, X= \sin \ov{\t} \d_{\t_1} - \cos \ov{\t} \d_{\t_2}, \d_{\ov{\t}} \}$. Let $J_{cyl}$ be the unique almost-complex structure satisfying $J(\d_t)= R_{\a_i} $ and $J(X)= \d_{\ov{\t}}$. One can readily check that $J_{cyl}$ is compatible with $d(e^t \a_i)$ and preserves $\op{ker}\a_i$.

Let $\ov{j}$ be the standard integrable complex structure on $S^2 \tms S^2$. Let $J$ be a compatible almost-complex structure on $(S^2 \tms S^2, \o_R \oplus \o_R)$ with the following two properties:
\begin{itemize}
\item[(i)] The restriction of $J$ to $\phi^i((-1,1) \tms S^*_{1, g_i}(\bb{T}^2))$ coincides with $(\phi^i)_*J_{cyl}$. 
\item[(ii)] J agrees with $\ov{j}$ in some open neighborhood $\mc{U}$ of $D_{\infty}$ which does not intersect $N(L_1) \cup N(L_2)$. 
\end{itemize}

We now introduce a family of compatible almost-complex structure $\{J_k^l\}$ on $(S^2 \tms S^2, \o_R \oplus \o_R)$ for $(k, l) \in \bb{N}_+ \tms \bb{N}_+$. We construct this family by stretching the neck along $S^*_{1, g_1}\bb{T}^2 \sub N(L_1)$ and $S^*_{1, g_2}\bb{T}^2 \sub N(L_2)$, following the procedure described in \cite[Sec. 2.7]{cmsftcompactness}. We will fix the convention that subscript indices correspond to stretching the neck along $S^*_{1,g_1}\bb{T}^2$ while the superscript indices correspond to stretching the neck along $S^*_{1,g_2}\bb{T}^2$. In other words, $J_k^l$ is obtained from $J$ by inserting a neck of length $k$ along $S^*_{1, g_1}\bb{T}^2$ and a neck of length $l$ along $S^*_{2, g_2}\bb{T}^2$. 

We also consider the almost-complex structures $\{J_{\infty}^l\}$, $\{J_k^{\infty}\}$ and $J_{\infty}^{\infty}$ on $S^2 \tms S^2-L_1$, on $S^2 \tms S^2 -L_2$ and on $S^2 \tms S^2-L_1-L_2$ respectively. These are constructed by replacing $J$ with $J_{cyl}$ in the image of $\phi^i((-\infty, 1) \tms S^*_{1, g_i}\bb{T}^2) \sub N(L_i) \sub S^2 \tms S^2-D_{\infty}$. 

\lem \label{lemma:perturbation} It is possible to choose $J$ such that the almost-complex structures $\{J_{\infty}^l \}, \{J_k^{\infty} \}$ and $J_{\infty}^{\infty}$ are regular for somewhere injective punctured curves. \elem

\pf Let $J_0$ be any compatible almost-complex structure on $S^2 \tms S^2$ which satisfies (i) and (ii) above. Observe that any compatible perturbation of $J_0$ which is fixed on $N(L_1) \cup N(L_2) \cup \mc{U}$ will also satisfy these properties. Let $\mc{V}= \mc{P}(R, R)- (N(L_1) \cup N(L_2) \cup \mc{U})$. 

Let $\{J(\a)\}_{\a=1}^{\infty}$ be an enumeration of the $\{ {J}_{\infty}^l \}, \{ {J}_k^{\infty} \}$ and ${J}_{\infty}^{\infty}$. 

It follows from elementary topological arguments that any punctured $J(\a)$-holomorphic curve which is somewhere injective must intersect $\mc{V}$. 

Observe that $J(\a)|_{\mc{V}}$ is independent of $\a$. For each $\a$, there is a Baire set of compatible perturbations of $J(\a)$ supported in $\mc{V}$, such that the perturbed almost-complex structure is regular for simply covered curves which intersect $\mc{V}$ (see \cite[Theorem 7.2]{wendlsft}). Since the $J(\a)$ are all equal inside $\mc{V}$, the space of perturbations of each $J(\a)$ can be naturally identified.

Since there are countably many $J(\a)$, the intersection of these Baire sets is nonempty. Hence there is a perturbation which works for all $\a$. If we apply this perturbation to $J_0$, then we obtain an almost-complex structure $J$ with the desired property. (Equivalently, we can think of this as simultaneously perturbing all of the $J(\a)$). \epf

Let $u$ be a punctured holomorphic curve mapping into $S^2 \tms S^2-L_1, S^2 \tms S^2-L_2$ or $S^2 \tms S^2-L_1-L_2$. We let $c_1^{\Phi}(u)$ be the relative Chern number of $u$ with respect to $\Phi= \{ \Phi^1, \Phi^2\}$. This is a count of zeros of a generic section of $u^*T(S^2 \tms S^2) \wedge u^*T(S^2 \tms S^2)$ which is constant near the punctures with respect to the trivializations induced by $\Phi^1$ and $\Phi^2$.

We have the following simple relation between the Chern number of a holomorphic plane and the Maslov index of its compactification.

\lem \label{lemma:maslovchern} For $i=1,2$, let $u_i: \C \to (S^2 \tms S^2 -L_i)$ be a $J$-holomorphic plane where $J= J_{\infty}^l$ or $J=J_k^{\infty}$. Let $v: \C \to S^2 \tms S^2 -L_1-L_2$ be a $J$-holomorphic plane for $J= J_{\infty}^{\infty}$. Then $2c_1^{\Phi}(u_i)= \mu(\ov{u_i})$ and $2c_1^{\Phi}(v) = \mu(\ov{v})$. \elem

\pf This is stated in \cite[Sec. 3.1 Eq. (2)]{igr} and references are provided for the proof. However, since these references follow notational and sign conventions which are different from ours, we will briefly sketch an argument in the appendix for the reader's convenience. \epf

\subsection{The moduli space of holomorphic planes} For a class $\a \in \pi_2(S^2 \tms S^2, L_1)$ and a distinguished point $p \in \C,$ we let 
\eqs \mc{M}_1(\a, J^l_{\infty}):= \{ u: \C \to S^2 \tms S^2-L_1 \mid \ov{\d}_{J^l_{\infty}} u =0, \ov{u}=\a \} / \{ \op{Aut}(\C, p) \} \eeqs
be the moduli space of $J^l_{\infty}$-holomorphic planes with one marked point whose compactification represents the class $\a.$ Let 
\begin{align*}
\op{ev}(\a, J_{\infty}^l): \mc{M}_1(\a, J_{\infty}^l) &\to S^2 \tms S^2-L_1 \\
([u], p) \mapsto u(p),
\end{align*} be the evaluation map. We will also denote by $\mc{M}(\a, J^l_{\infty})$ the moduli space of $J^l_{\infty}$-holomorphic planes (with no marked points) whose compactification represents the class $\a.$

If we assume that $\a \in \pi_2(S^2 \times S^2, L_1)$ is primitive, it follows from \Cref{lemma:perturbation} that $\mc{M}_1(\a, J_{\infty}^l)$ and $\mc{M}(\a, J_{\infty}^l)$ are a smooth manifolds and that $\op{ev}(\a, J_{\infty}^l)$ is a smooth map. 

We now come to the following important proposition. 

\prop \label{proposition: igrciting} For every natural number $l>0,$ there exists a class $\a_l \in \pi_2(S^2 \tms S^2, L_1)$ with $\mu(\a_l)=2$ such that the following properties are satisfied:
\begin{itemize}
\item[(i)] There exists a component $\mc{M}^0(l) \sub \mc{M}(\a_l, J_{\infty}^l)$ diffeomorphic to $S^1,$ and a diffeomorphism $\mc{M}^0_1(\a_l, J^l_{\infty}) \simeq \mc{M}^0(l) \tms \C.$ Here $\mc{M}_1^0(\a_l, J_{\infty}^l)$ is a component of $\mc{M}_1(\a_l, J_{\infty}^l),$ the moduli space of planes in the class $\a_k$ with one marked point. 
\item[(ii)] The evaluation map $\op{ev}(\a_l, J_{\infty}^l): \mc{M}^0_1(l) \tms \C \to S^2 \tms S^2 -L_1$ is a smooth embedding, and its image is disjoint from $D_{\infty}$. Thus, we can also view $\op{ev}(\a_l, J_{\infty}^l)$ as mapping into $\mc{P}(R, R) -L_1 \sub \R^4-L_1$. 
\item[(iii)] The evaluation map can be modified to a yield smooth map $\ov{\op{ev}}(\a_l, J_{\infty}^l): (\mc{M}^0(l) \tms \ov{D^2}, \mc{M}^0(l) \tms S^1) \to (\R^4, L_1)$, whose image can be made to lie in an arbitrarily small neighborhood of the image of $\op{ev}(\a_l, J_{\infty}^l)$.
\end{itemize}
The planes belonging to the component $\mc{M}^0(l)$ will be called ``small planes"; cf. \cite[Sec. 5.2]{igr}. 
\eprop

\rmk The choice of the class $\a_l$ and component $\mc{M}^0(l)$ in \Cref{proposition: igrciting} is not canonical. In general, there could be multiple families of small planes. \ermk

\pf This proposition follows from the analysis carried out in \cite{igr}. It follows from a well-known theorem of Gromov that $S^2 \tms S^2$ is foliated by ${J}_k^l$-holomorphic spheres in the classes $[S^2 \tms *]$ and $[* \tms S^2]$. If one views $l$ as fixed and sends $k \to \infty$, one can analyze the limiting behavior of these spheres under the SFT compactness theorem. This analysis is carried out in \cite[Sec. 5]{igr}. We observe in light of \Cref{lemma:perturbation} that the almost-complex structures $\{{J}_{\infty}^l\}$ satisfy the transversality properties which are assumed in this analysis (see \cite[p. 27]{igr}).

It follows from \cite[Prop. 5.11]{igr} that there is a component $\mc{M}^0(l) \sub \mc{M}(\a_l, J_{\infty}^l)$ satisfying (i) and (ii) of \Cref{proposition: igrciting}. These planes are referred to in \cite{igr} as ``small planes", and we will continue to use this terminology.  

The proof of (iii) is carried out in \cite[Sec. 5.3]{igr}. The key input is \cite[Lem. 5.13]{igr} which guarantees that distinct planes are asymptotic to distinct Reeb orbits. The desired modification can then be constructed using a standard asymptotic formula for punctured holomorphic curves, as in \cite[Sec.\ 5.3]{igr}. \epf

We record the following lemma which will be useful in the next section. 

\lem \label{lemma:awaydivisor} There exists a neighborhood $\mc{U}' \sub \mc{U}$ independent of $l \in \bb{N}_+$ with $D_{\infty} \sub \mc{U}'$ such that none of the $J_{\infty}^l$-holomorphic small planes intersect $\mc{U}'$. \elem
\pf Let $O \sub S^2$ be a small open neighborhood containing $\{\infty\} \in S^2$, with the property that $\{p\} \tms S^2 \sub \mc{U}$ and $S^2 \tms \{ p' \} \sub \mc{U}$ for all $p, p' \in O$. Since the asymptotic boundaries of the small planes are geodesics of $L_1$ and since $N(L_1) \cap \mc{U}= \emptyset$, it follows that the intersection number of a small plane with the spheres $\{p\} \tms S^2$ and $S^2 \tms \{ p' \}$ is independent of $p, p' \in O$. This intersection number must be zero since the small planes do not intersect $D_{\infty}$. Since $J$ is standard in $\mc{U}$, it follows by positivity of intersection that the small planes do not intersect $\{p\} \tms S^2$ and $S^2 \tms \{ p' \}$ for any $p, p' \in O$.\epf

\subsection{Deforming the complex structure by stretching the neck} We now implement the second part of the proof of \Cref{theorem: erimprove}. We will show that the moduli spaces considered in \Cref{proposition: igrciting} eventually become disjoint from $L_2$.

\prop \label{proposition: sim1} There exists $\L \gg 1$ such that the image of $\op{ev}(\a_l, J_{\infty}^l): \mc{M}^0(l) \tms \C \to S^2 \tms S^2 -L_1$ is disjoint from $L_2$ for all $l>\L.$\eprop

\pf The proof is similar to \cite[Theorem 4.1]{evans-rizell}. Let us suppose for contradiction that the statement is false. Then there exists a sequence $\{ u_l \}$ of $J_{\infty}^l$ holomorphic planes in the class $\a_l$ such that $\op{Im} u_l \cap L_2$ is nonempty for all $l.$  Recall by \Cref{proposition: igrciting} (ii) that $u_l \cap D_{\infty} = \emptyset$, so we can view $\op{ev}(\a_l, J_{\infty}^l)$ as mapping into $\mc{P}(R, R) -L_1 \sub \R^4-L_1$. Since $L_1$ is monotone as a Lagrangian submanifold of $\R^4$, there is a constant $C>0$ such that $\o(u_l)= C\mu(\a_l)= 2C.$ Up to replacing $\{u_l\}$ with a subsequence, it follows by the by the SFT compactness theorem (cf. \Cref{remark:sft}) that the sequence $\{u_l\}$ converges to a holomorphic building $\mathbf{u}.$

For $\s=1,2,\dots,N,$ let $\{u^{\s}\}$ be an enumeration of the components of $\mathbf{u}$. The $u^{\s}$ map into domains which are diffeomorphic to $S^2 \tms S^2-L_1-L_2, (\R \tms S^*_{1,g_1}\bb{T}^2), (\R \tms S^*_{1,g_2} \bb{T}^2)$ and $T^*L_2.$ In light of \Cref{lemma:awaydivisor}, the planes $u_l$ stay uniformly away from $\d \ov{\mc{P}(R,R)}$. This implies that the $u^{\s}$ which map into $S^2 \tms S^2 - L_1-L_2$ actually land inside $\mc{P}(R,R)-L_1-L_2 \sub \R^4-L_1-L_2$.

Since $\R^4-L_1-L_2$ is exact, it follows that any $u^{\s}$ mapping into $\R^4-L_1-L_2$ has at least one puncture. The domains $(\R \tms S^*_{1,g_1}\bb{T}^2), (\R \tms S^*_{1,g_2} \bb{T}^2)$ and $T^* L_2$ have vanishing homotopy groups in degree strictly greater than $1,$ so all $u^{\s}$ mapping into these domains must also have at least one puncture.

It now follows by elementary topological considerations that the building $\mathbf{u}$ must contain a plane. Up to relabeling the indices, we can assume that $u^1$ is a plane. Observe that $u^1$ must map into $\R^4-L_1-L_2$ due to the fact that the flat metric $g_i$ on $L_i$ admits no contractible geodesics for $i=1,2$. 

Let $\ov{u}^1$ be the compactification of $u^1.$ By combining \Cref{lemma: next} and \Cref{lemma:posindex} below, we find that $\ov{u}^1$ has Maslov index $2$. Hence $u^1$ cannot converge at its puncture to a geodesic of $(L_2, g_2)$ since $\o (u_k)< A_2(L_2).$ Hence it converges at its puncture to a geodesic of $(L_1, g_1)$. This implies that $\ov{u^1}$ has area $A_2(L_1)=\o(u_k)$ since $L_1$ is monotone. It follows that there are in fact no other components to the building $\bf{u},$ which contradicts the assumption that the $u_l$ intersect $L_2.$ \epf

\rmk  \label{remark:sft} In the proof of \Cref{proposition: sim1}, we are appealing to a version of the SFT compactness theorem for ``neck-stretching" in a manifold with a negative cylindrical end. To the author's knowledge, a proof of this precise version of the SFT compactness theorem does not appear in the literature, but closely related versions are described in \cite{sftcompactness} and the arguments there go through in our setting with straightforward modifications. We note that an alternative approach to SFT compactness is detailed in \cite{cmsftcompactness}. 

It seems useful to clarify the relation between the symplectic area of the holomorphic planes $u_l$ considered in the proof of \Cref{proposition: sim1}, and the notions of energy considered in \cite{sftcompactness} which are needed for proving compactness. Although these notions are strictly different, it can be shown using the arguments of \cite[Lem.\ 9.2]{sftcompactness} that the symplectic area of the planes $u_l$ controls the relevant energies in \cite{sftcompactness}. For completeness, the details of this argument are provided in the Appendix; see \Cref{subsection:energyappendix}. \ermk 

\rmk If we assume that $L_2$ is monotone, it follows that the disk $\ov{u}_1$ considered in the proof of \Cref{proposition: sim1} above has Maslov index at least $2$. This makes it possible to prove \Cref{proposition: sim1} without appealing to \Cref{lemma: next} and \Cref{lemma:posindex}. The reader who is only interested in monotone tori (and in particular in the proof of \Cref{theorem:clifforddisjoint}) can therefore safely pass to \Cref{section: enumerativeinvariant} of this paper. \ermk

\lem[cf. Prop. 3.5 in \cite{igr}] \label{lemma: next} The sum of the Fredholm indices of the components of the building $\mathbf{u}$ which map into $\R^4-L_1-L_2$ is at most $1.$ \elem

\pf As in the proof of \Cref{proposition: sim1}, let $\{u^{\s}\}$ be an enumeration of the components of the building $\mathbf{u}$ for $\s=1,2,\dots,N$. Let $K$ be the total number of asymptotic Reeb orbits of the components of the building $\mathbf{u}.$ Let us compute the sum of the Fredholm indices of all the $u^{\s}.$ We claim that the following equation holds
\eq \label{equation: all-levelsempty} \sum_{\s=1}^N \op{ind}(u^{\s}) = -2N + (3K-2) + 2\sum_{\s=1}^N c_1^{\Phi}(u^{\s}) = -2N+(3K-2) + 2c_1^{\Phi}(u_l), \eeq where we assume that $l$ is large enough so that $c_1^{\Phi}(u_l)$ is independent of $l$. 

The key input in proving \eqref{equation: all-levelsempty} is the index formula \eqref{equation: explicitindex} in the appendix. Observe that, with a single exception, all asymptotic orbits of the components of $\mathbf{u}$ occur as the negative puncture of exactly one component and as the positive puncture of exactly one component. Let $\rho$ be the unique asymptotic orbit which is not paired up with a positive puncture. It follows that every asymptotic orbit gets counted three times in \eqref{equation: explicitindex}, except for $\rho$ which gets counted only once. 

In fact, since $\mathbf{u}$ is a limit of planes, it follows that $K=N$ and hence 
\eq \label{equation: all-levels1empty} \sum_{\s=1}^N \op{ind}(u^{\s})  = N-2+ 2 c_1^{\Phi}(u_l). \eeq

We now argue as in \cite[Lemma 3.1]{igr}. Let $T \sub \{1,2,\dots,N\}$ be the set of all $\s \in \{1,2,\dots,N\}$ such that $u^{\s}$ maps into the domains $T^*L_2, \R \tms S^*_{1, g_1} \bb{T}^2$ or $S^*_{1, g_2} \bb{T}^2$. With the exception of $\rho$, all asymptotic orbits of the components of $\mathbf{u}$ occur as a positive puncture of some $u^{\s}$ for $\s \in T$. Hence it follows from \eqref{equation: explicitindex} that that the set of punctured curves $\{u^{\s}\}_{\s \in T}$ has total index \eq \label{equation: lower-levelsempty} \sum_{\s \in T} \op{ind}(u^{\s}) = - \sum_{\s \in T} \chi(u^{\s})+ (K-1) = -\sum_{\s \in T} \chi(u^{\s}) + (N-1). \eeq 

Combining \eqref{equation: all-levels1empty} with \eqref{equation: lower-levelsempty}, we find that the sum of the indices of the components mapping into $\R^4-L_1-L_2$ is precisely
\begin{align} \sum_{\s=1}^N \op{ind}(u^{\s}) - \sum_{\s \in T} \op{ind}(u^{\s}) &= N-2 + 2c_1^{\Phi}(u_l)+ \sum_{\s \in T} \chi(u^{\s}) -(N-1) \\
&\leq -1+ 2c_1^{\Phi}(u_l) = -1+ \mu(\ov{u_l})=1, \nonumber
\end{align}
where we have used \Cref{lemma:maslovchern} and the previously observed fact that all components mapping into $T^* \bb{T}^2$, $\R \tms S^*_{1,g_1} \bb{T}^2$ or $\R \tms S^*_{1,g_2} \bb{T}^2$ must have at least two punctures and thus have non-positive Euler characteristic.  \epf

\lem[cf. Lem. 3.1 and Lem. 3.3 in \cite{igr}] \label{lemma:posindex} Suppose that $u^{\tau} \in \{u^{\s}\}_{\s=1}^N$ is a component of the building $\mathbf{u}$. Then $\op{ind}(u^{\tau}) \geq 0$. If $u^{\tau}$ is a plane, then $\op{ind}(u^{\tau}) \geq 1$ with equality if and only if $u^{\tau}$ is simply-covered and the compactification $\ov{u}^{\tau}$ has Maslov index $2$. \elem

\pf If $u^{\tau}$ maps into $\R \tms S^*_{1,g_i} \bb{T}^2$ for $i=1,2$, then $\chi(u^{\tau}) \leq 0$ (since we saw that $u^{\tau}$ has at least 2 punctures) and $c_1^{\Phi}(u^{\tau})=0$. It follows that $\op{ind}(u^{\tau}) = -\chi(u^{\tau})+ 2c_1^{\Phi}(u^{\tau}) \geq 0$. 

We can therefore assume that $u^{\tau}$ maps into $\R^4-L_1-L_2$. In this case, recall from the proof of \Cref{proposition: sim1} that $u^{\tau}$ has at least one puncture. It follows from \Cref{lemma:perturbation} that $\op{ind}(u^{\tau}) \geq 0$ if $u^{\tau}$ is simply-covered.  

If $u^{\tau}$ is multiply covered, then there exists a map $\phi: \dot{\S} \to \dot{\S}'$ such that $\op{deg}(\phi) = d>1$ and $u^{\tau} = v^{\tau} \circ \phi$ where $v^{\tau}$ is a simply-covered punctured curve. Here $\dot{\S}$ and $\dot{\S}'$ are punctured Riemann surfaces of genus $0$, having $k_{u^{\tau}}$ and $k_{v^{\tau}}$ punctures respectively.

Let $\mc{B}$ be an enumeration of the branch points of $\phi$ and set \eqs b= \sum_{p \in \mc{B}} (m_p-1),\eeqs where $m_p$ is the multiplicity of $\phi$ at $p$. By the Riemann-Hurwicz formula, we have \eq 2=2d-b.\eeq Observe that we also have \eq d k_{v^{\tau}} \leq k_{u^{\tau}} + b. \eeq

In light of the index formula \eqref{equation: explicitindex}, we now have 
\begin{align*}
\op{ind}(u^{\tau}) &= -2 + k_{u^{\tau}} + 2 c_1^{\Phi}(u^{\tau}) \\
&= b-2d+k_{u^{\tau}}+  2 c_1^{\Phi}(u^{\tau})  \\
& \geq d k_{v^{\tau}} -2d + 2d c_1^{\Phi}(v^{\tau}) \\
&= d \op{ind}(u^{\tau}).
\end{align*}

Assuming that $\op{ind}(u^{\tau}) <0$, it then follows that $\op{ind}(v^{\tau}) <0$. This is a contradiction since $v^{\tau}$ is simply-covered. This proves the first part of the lemma.  

If $u^{\tau}$ is a plane, then it follows from \Cref{lemma:maslovchern} that $\op{ind}(u^{\tau})= -1+ \mu(\ov{u}^{\tau}) \geq 1$, with equality if and only if $\mu(\ov{u}^{\tau})=2$. \epf

\pf[Proof of \Cref{theorem: erimprove}] The theorem follows immediately by combining \Cref{proposition: sim1} and \Cref{proposition: igrciting} (iii). \epf

\section{Enumerative invariants and linking obstructions} \label{section: enumerativeinvariant}

Let $L_1$ and $L_2$ be monotone Lagrangian tori and assume that the monotonicity factor of $L_1$ is at most equal to that of $L_2$. It follows from \Cref{corollary: pi1linking} that the only obstruction to smoothly unlinking $L_1$ and $L_2$ is the nontriviality of the map $\pi_1(L_1) \to \pi_1(\R^4-L_2)$. In this section, we will relate this map to certain enumerative counts of holomorphic disks with boundary in $L_1$. 

\subsection{An enumerative invariant of Lagrangian tori}  \label{subsection: enumerativeinvariant1} The present section follows \cite[Section 3]{auroux}. Fix a monotone Lagrangian torus $L \sub (\R^4, \o)$ and let $J$ be a compatible almost-complex structure which is standard at infinity.  Throughout this section, all almost-complex structures (and families of almost-complex structures) will be assumed to coincide with the standard complex structure $j$ at infinity. 

Given a class $\a \in \pi_2(\R^4, L)$, let \eqs \mc{M}(\a, J)= \{ u: (D^2, \d D^2) \to (\R^4, L) \mid \ov{\d}_J u =0, u_*[D^2]= \a \} \eeqs be the moduli space of $J$-holomorphic disks representing $\a$. Let \eqs \mc{M}_{0,1}(\a, J)= \{ u: (D^2, \d D^2) \to (\R^4, L) \mid \ov{\d}_J u =0, u_*[D^2]= \a \} / \op{Aut}(D^2,1)\eeqs be the moduli space of $J$-holomorphic disks representing $\a$ with one boundary marked point.

If $\a$ is primitive and $J$ is regular for simply-covered curves, then $\mc{M}(\a, J)$ and $\mc{M}_{0,1}(\a, J)$ are smooth manifolds of dimension $-1+ \mu(\a)$ and $ \mu(\a)$ respectively.  The boundary evaluation map 
\begin{align*}
\op{ev}(\a, J) : \mc{M}_{0,1}(\a, J)  &\to L \\
[u] &\mapsto u(1)
\end{align*} is also smooth. 

\defi \label{definition:enumerativeinvariant} Let $\a \in \pi_2(\R^4,L)$ be a class with $\mu(\a)=2$ and let $J$ be a compatible almost-complex structure which is regular for simply-covered curves with boundary in $L$. We define $n(L, \a) \in \Z/2$ to be the mod $2$ degree of the boundary evaluation map $\op{ev}(\a, J): \mc{M}_{0,1}(\a, J) \to L.$ \edefi

The invariant $n(L, \a)$ can be interpreted as a count of holomorphic disks representing the class $\a$ which pass through a generic point of $L.$ A standard cobordism argument (which uses crucially our assumption that $L$ is monotone) shows that $n(L,\a)$ remains unchanged under Hamiltonian isotopies of $L,$ and under generic homotopies between regular almost-complex structures. Since any two regular almost-complex structures can be connected by a generic homotopy, it follows that $n(L, \a)$ is independent of the choice of regular almost-complex structure. 

\ex \label{example:cliffordenumerative} Consider the Clifford torus $L= L(r,r) = \{ (z_1, z_2) \in \C^2 \mid |z_1|=|z_2|=r\}$ for some $r>0$. It can be shown (cf. \cite[Thm. 10.2]{cho}) that the standard complex structure $j$ is regular for all holomorphic disks with boundary in $L$. Let $\a_1= [D^2 \tms *]$ and let $\a_2= [* \tms D^2]$ be classes in $\pi_2(\R^4, L)$. It's clear that $\op{ev}(\a_i, j)$ is a degree $1$ map for $i=1,2$.  It follows that $n(L, \a_i)=1$. \eex

\subsection{Application of the invariant to linking} \label{subsection:applicationtolinking} In this section, we explain why the enumerative invariant introduced above is relevant for the study of linking of Lagrangian tori. In particular, we will use it to show that any two Clifford tori in $(\R^4, \o)$ are always smoothly unlinked, thus proving \Cref{theorem:clifforddisjoint} in the introduction. 

We begin with the following key proposition. 

\prop \label{proposition: essential} Let $L_1 \sub (\R^4, \o)$ be a monotone Lagrangian torus and let $L_2 \sub (\R^4, \o)$ be a (not necessarily monotone) Lagrangian torus disjoint from $L_1.$  Suppose that there exist homotopy classes $\a_1, \a_2 \in \pi_2(\R^4,L_1)$ satisfying the following properties:
\begin{enumerate}
\item[(i)] $\mu(\a_1)=\mu(\a_2)=2,$
\item[(ii)] $ n(L_1, \a_1)= n(L_1, \a_2) = 1, $
\item[(iii)] The image of $\{\a_1, \a_2\}$ under the inclusion $\pi_2(\R^4, L_1) \xrightarrow{\sim} \pi_1(L_1) \to \pi_1(L_1) \otimes \Q$ generates a basis.
\end{enumerate}

If $A_2(L_2) \geq A_2(L_1),$ then the group homomorphism $i_*: \pi_1(L_1) \to \pi_1(\R^4-L_2)$ induced by the inclusion $i: L_1 \to \R^4-L_2$ is trivial. \eprop

\pf  It follows from \Cref{lemma: homotopyofcomplement} that $\pi_1(\R^4-L_2)$ is torsion-free. Hence, in light of property (iii) above, it is enough to prove that $\a_1$ and $\a_2$ have trivial image in $\pi_1(\R^4-L_2).$  

Fix a compatible almost-complex structure $J$ on $(\R^4, \o)$. As in \Cref{subsection:neckstretch}, let $\{J^l\}_{l=1}^{\infty}$ be a sequence of almost-complex structures obtained from $J$ by stretching the neck in a Weinstein neighborhood $N(L_2)$ which is disjoint from $L_1.$  Let $N(L_1)$ be a Weinstein neighborhood of $L_1$ with the property that $N(L_1) \cap N(L_2)= \emptyset.$  

It follows from property (i) and from the fact that $L_1$ is orientable that $\a_1$ and $\a_2$ are primitive classes. By standard genericity arguments, we can perturb $J$ outside of $N(L_1) \cup N(L_2)$ in such a way that all $J^l$ can be assumed to be regular for simply-covered pseudoholomorphic disks; cf. \Cref{lemma:perturbation}.

Since $n(L_1, \a_i)=1$ for $i=1,2,$ it follows that $\mc{M}(\a_i, J^l)$ is non-empty for all $l \geq 1.$  The proposition is now a consequence of the following lemma.  \epf

\lem \label{lemma:sim11} For $i=1,2,$ there exists a large integer $\Lambda$ such that any disk $u_l \in \mc{M}(\a_i, J^l)$ is disjoint from $L_2$ if $l>\Lambda.$ \elem

\pf The proof is similar to that of \Cref{proposition: sim1} and \cite[Theorem 4.1]{evans-rizell}. Suppose for contradiction that the statement is false. This implies that there exists an infinite sequence of $J^l$-holomorphic disks $u_l$ such that $u_l \cap L_2$ is non-empty for all $l \in \N_+.$ Up to passing to a subsequence, we can assume by the SFT compactness theorem that the $u_l$ converge to a holomorphic building $\mathbf{u}$. This building must have a component mapping into the domain $T^*L_2$ due to our assumption that $u_l \cap L_2$ is non-empty. 

It can be shown by a routine modification of the proofs of \Cref{lemma: next} and \Cref{lemma:posindex} that the components of $\mathbf{u}$ satisfy the following two properties.
\begin{itemize}
\item[(i)] The sum of the Fredholm indices of the components of $\mathbf{u}$ which map into $\R^4-L_2$ is at most $1$. 
\item[(ii)] Every component $u^{\tau}$ of $\mathbf{u}$ has non-negative Fredholm index. If $u^{\tau}$ is a plane, then $\op{ind}(u^{\tau}) \geq 1$ with equality if and only if $u^{\tau}$ is simply-covered and the compactification $\ov{u}^{\tau}$ has Maslov index $2$. 
\end{itemize}

It now follows by an argument analogous to the proof of \Cref{proposition: sim1} that the limit building $\mathbf{u}$ must contain a plane $v$ whose compactification $\ov{v}$ has Maslov index $2.$ Observe that $\ov{v}$ cannot have boundary in $L_2$.  Indeed, we must have $\omega(\ov{v}) \geq A_2(L_2) \geq A_2(L_1)= \omega(u_l)$, which would imply that there are no other components of the holomorphic building mapping into $\R^4-L_2$. But there must be at least one other component of the building having boundary in $L_1$.  It follows that $\ov{v}$ has boundary in $L_1$. But this means that $\omega(\ov{v}) \geq A_2(L_1)= \omega(u_l)$. Hence there are no other components of the building mapping into $\R^4-L_2$. By connectedness of the building (cf. Sec.\ 9.1 (v) of \cite{sftcompactness}), it follows that there are no components of the building mapping into $T^*L_2$, which is a contradiction in view of the paragraph above. \epf

\rmk Note that for the purpose of proving \Cref{theorem:clifforddisjoint} in the introduction, we could bypass the above argument entirely by assuming $L_2$ to be monotone. In this case, the proof is identical to the original argument of Dimitroglou Rizell and Evans \cite[Thm. 4.1]{evans-rizell}. \ermk

We arrive at the following corollaries, the second of which implies \Cref{theorem:clifforddisjoint} in the introduction. 

\cor \label{corollary:cliffordunlink2} Suppose that $L_1, L_2 \sub (\R^4, \o)$ are disjoint Lagrangian tori satisfying the assumptions of \Cref{proposition: essential}. Suppose also that $A_2(L_2) \geq A_2(L_1)$. Then $L_1$ and $L_2$ are smoothly unlinked. \ecor
\pf It follows from \Cref{proposition: essential} that $\pi_1(L_1)$ has trivial image in $\pi_1(\R^4-L_2).$ By \Cref{corollary: pi1linking}, this implies that $L_1$ and $L_2$ are smoothly unlinked. \epf

\cor \label{corollary:cliffordunlink1} Any two Clifford tori (of possibly different monotonicity factor) in $\R^4$ are smoothly unlinked. \ecor
\pf If $L_1$ and $L_2$ are both Clifford tori, we can assume without loss of generality that $A_2(L_2) \geq A_2(L_1)$. It follows from \Cref{example:cliffordenumerative} that $L_1$ satisfies the assumptions of \Cref{proposition: essential}. Hence \Cref{corollary:cliffordunlink1} follows from \Cref{corollary:cliffordunlink2}. \epf

\subsection{Configurations of monotone Lagrangian tori} 

The goal of this section is to characterize possible configurations of monotone Lagrangian tori in $\R^4$ up to smooth isotopy. Some of the arguments will only be sketched, since they are not needed in the remainder of this paper. 

We begin with an auxiliary lemma. 

\lem Let $\l= x_1 dy_1+x_2 dy_2.$ Let $\g_1$ and $\g_2$ be simple closed curves in $\R^4.$ Then $\g_1$ and $\g_2$ are Hamiltonian isotopic if and only if $\int_{\g_1} \l = \int_{\g_2} \l.$ \elem

\pf The main step is to observe that there exists a smooth isotopy $\{\g_t\}$ such that $\int_{\g_t} \l$ is independent of $t.$  By the symplectic neighborhood theorem, we can then extend the isotopy to a compactly supported diffeomorphism $\Phi$ which is a symplectomorphism near $\g_t.$  Let $\o_t= \Phi^* \o$ and observe that it would be enough to produce a compactly supported isotopy $\Psi_t$ such that $\Psi_t^* \o_t= \o.$  This can be accomplished by a standard Moser-type argument, which relies on the fact that the $\g_t$ have constant action. \epf

\cor \label{corollary: chekanovsqueeze} Let $\g \sub \R^4$ be a simple closed curve and let $\g \sub \mc{U}$ be a tubular neighborhood.  Then there exists a Chekanov torus $L_{Ch} \sub \mc{U}$ such that the map $\pi_1(L_{Ch}) \to \pi_1(\mc{U})=\Z$ is surjective. \ecor
\pf Choose a simple closed curve $\tilde{\g} \sub \mc{U}$ such that $[\g]= [\tilde{\g}] \in \pi_1(\mc{U})$ and $\int_{\tilde{\g}} \lambda =0$. By the lemma, there exists a global Hamiltonian isotopy taking the curve $\kappa(t)= (\cos t, 0, \sin t, 0)$ to the curve $\tilde{\g}.$ In particular, this isotopy maps a small neighborhood $\mc{V}$ of $\kappa$ into $\mc{U}$. Inspecting the definition of the Chekanov torus in \Cref{example: chekanov}, we can choose the monotonicity factor small enough so that $L_{Ch}(r^2) \sub \mc{V}.$ The corollary follows.

The existence of $\tilde{\g}$ is geometrically clear but tedious to prove in detail. A sketch of a possible argument goes as follows. Pick a point $p \in \g$ and a small ball $B_p \sub \mc{V}$. Supposing that $\int_{\g} \l=A$, one can clearly construct an immersed closed curve $c_p: [0,1] \to B_p \sub \R^4$ with $c_p(0)=c_p(1)=p$ and such that $\int_{c_p} \l= -A$. The concatenation $\g * c_p$ is now an immersed closed curve of area zero satisfying the desired properties. By wiggling it slightly, we can get a nearby embedded curve of area zero.  \epf

\ex \label{example:chekanovlink} Let $L_1 \sub (\R^4, \o)$ be an arbitrary Lagrangian torus. Let $\g \sub \R^4-L_1$ be a simple closed curve which realizes a nontrivial element of $\pi_1(\R^4-L_1) \simeq \Z$; cf. \Cref{lemma: homotopyofcomplement}. Let $\mc{U} \sub \R^4-L_1$ be a tubular neighborhood of $\g$. It follows from \Cref{corollary: chekanovsqueeze} that there is a Chekanov torus $L_2 \sub \mc{U} \sub \R^4-L_1$ with the property that $\pi_1(L_2) \to \pi_1(\R^4-L_1)$ has nontrivial image. By \Cref{lemma: pi1link}, $L_1$ is homologically linked with $L_2$. In particular, $L_1$ and $L_2$ are not smoothly unlinked; cf. \Cref{corollary:cliffordunlink1}. \eex

Let $\mc{C}$ be a finite collection of disjoint monotone Lagrangian tori in $\R^4$. Let $a_1>\dots>a_n$ be the set of values of $A_2(L)$ for $L \in \mc{C}$. We partition $\mc{C}$ into \emph{levels} $\ell_1,\dots,\ell_n$ by stipulating that $L \in \ell_i$ if  $A_2(L)= a_i$. 

Observe that $\mc{C}$ satisfies the following properties:
\begin{itemize}
\item[(i)] All pairs of Clifford tori in $\mc{C}$ are smoothly unlinked from one another. This follows from \Cref{corollary:cliffordunlink1}. 
\item[(ii)] All tori of the same level are smoothly unlinked. Any torus $L \in \mc{C}$ bounds a solid torus embedded in the complement of all tori of higher level. This follows from \Cref{corollary:restateB}.
\end{itemize}

We can think of $\mc{C}$ as being built in the following way. First, one chooses disjoint monotone tori $L_1^1,\dots,L_1^{j_1}$ with $A_2(L_1^i)=a_1$ which form the level $\ell_{1}$. Having constructed the levels $\ell_1,\dots,\ell_{n-1}$, we construct the level $\ell_n$ by choosing disjoint monotone tori $L_n^1,\dots,L_n^{j_n}$ such that $A_2(L_n^i)= a_n$. We require that the $L_n^k$ do not intersect any of the previously constructed $\ell_i$.

At each step, one can consider the image of the map $\pi_1(L_n^i) \to \pi_1(\R^4- \cup_{1}^{n-1} \ell_k)$, which is a cyclic subgroup. These subgroups are discrete invariants of our construction of $\mc{C}$. The following proposition shows that they are in a sense the only invariants of the construction.
\prop[Uniqueness] For $i=1,2,\dots,j_n$, the Lagrangian torus $L_n^i$ is entirely determined up to smooth isotopy by the image of the map $\pi_1(L_n^i) \to \pi_1(\R^4- \cup_{1}^{n-1} \ell_k)$. \eprop
\pf The arguments of \Cref{subsection:neckstretch} allow us to produce a smoothly embedded solid torus which does not intersect any of the tori belonging to the levels $\ell_1,\dots,\ell_{n-1}$. This can be done by stretching the neck along all of these tori simultaneously. The remainder of the proof is now analogous to the proof of \Cref{corollary: pi1linking}. \epf

The next proposition shows that all possible images of the maps $\pi_1(L_n^i) \to \pi_1(\R^4- \cup_{1}^{n-1} \ell_k)$ are indeed achieved through the above construction. 
\prop[Existence] Let $S \sub \pi_1(\R^4- \cup_1^{n-1} \ell_k)$ be a cyclic subgroup. Then there exists a torus $L$ such that $\pi_1(L) \to \pi_1(\R^4- \cup_1^{n-1} \ell_k)$ has image $S$ and $A_2(L)< a_i$ for $i=1,\dots,n-1$. \eprop

\pf Let $\g \sub \R^4-\cup_1^{n-1} \ell_k$ be a simple closed curve generating $S$. By modifying $\g$ in a $C^0$-small neighborhood, we can assume that $\int_{\g} \l=0$.  Let $\mc{U} \sub \R^4-\cup_1^{n-1} \ell_k$ be a tubular neighborhood of $\g$. It now follows from \Cref{corollary: chekanovsqueeze} that there is a Chekanov torus $L \sub \mc{U}$ such that $\pi_1(L) \to \pi_1(\R^4- \cup_1^{n-1} \ell_k)$ has image precisely $S$. \epf

The upshot of the above propositions is that monotone Lagrangian tori in $\R^4$ are essentially characterized up to smooth isotopy by a discrete set of topological choices.  In fact, by a repeated application of the arguments of \Cref{corollary: pi1linking}, one should be able to prove a statement to the effect that isomorphic choices of this data give rise to smoothly isotopic configurations of tori. We leave it to the interested reader to formulate a precise version of this statement.

\section{Homological linking of non-monotone tori in $\R^4$} \label{section: non-monotone}

In this section, we introduce a class of non-monotone Lagrangian tori in $\R^4$ whose members will be called \emph{admissible tori}. We show that this class is closed under Hamiltonian isotopies and contains ``most" product tori. The main result of this section is \Cref{theorem: nonmonotone2} (stated as \Cref{theorem: nonmonotone} in the introduction), which gives sufficient conditions under which admissible tori are homologically unlinked and thus modestly generalizes \Cref{theorem:edr} of Dimitroglou Rizell and Evans in dimension $4$. We will show in \Cref{section: linkingconstruction} that \Cref{theorem: nonmonotone2} is sharp in an appropriate sense.

\subsection{An enumerative invariant for admissible Lagrangian tori} Let $L \sub (\R^4, \o)$ be a Lagrangian torus. We will assume throughout this section that $L$ is not monotone; this is a harmless assumption since the results proved in this section will be weaker than those of the previous two sections, which do apply to monotone tori.

Unless otherwise indicated, all almost-complex structures in this section are assumed to coincide at infinity with the standard complex structure $j$. Let $J$ be an almost-complex structure on $\R^4$ which is compatible with $\o$ and regular for simply covered curves.  Given a primitive class $\a \in \pi_2(\R^4-L),$ let $\mc{M}(\a, J), \mc{M}_{0,1}(\a, J)$ and $\op{ev}(\a, J)$ be defined as in \Cref{section: enumerativeinvariant}.

Observe that there is a unique class $\a_0 \in \pi_2(\R^4, L)$ with the property that $\mu(\a_0)=2$ and $\o(\a_0)=A_2(L)$. The existence of this class follows from the fact that every Lagrangian torus in $\R^{2n}$ admits a disk of Maslov index $2$ (this was proved by Cieliebak and Mohnke \cite[Theorem 1.2]{cmaudin}, although the $4$-dimensional case was already known). The uniqueness of this class follows from our assumption that $L$ is not monotone. 

\defi \label{definition:muinfimal} Let $L \sub (\R^4, \o)$ be a Lagrangian torus and let $\a_0 \in \pi_2(\R^4, L)$ be the unique class with the property that $\mu(\a_0)=2$ and $\o(\a_0)=A_2(L)$. By analogy with \cite[Def. 4.1]{evans-rizell}, we will call $\a_0$ the \emph{$\mu$-infimal} class of $L$. \edefi

\defi Let $\a \in \pi_2(\R^4, L)$ be a primitive class and let $J$ be an $\o$-compatible almost complex structure which is regular for simply-covered curves. We define $\tilde{n}(L, \a, J) \in \Z/2$ to be the mod $2$ degree of the evaluation map $\op{ev}(\a, J): \mc{M}_{0,1}(\a, J) \to L$. \edefi

Since $L$ is not monotone, one does not in general expect the count $\tilde{n}(L, \a, J)$ to be independent of $J$. However, we have the following useful proposition.

\prop \label{proposition:Jindep} Let $\a_0$ be the $\mu$-infimal class of $L$. Then $\tilde{n}(L, J, \a_0) \in \Z/2$ is independent of the choice of $J$ among $\o$-compatible almost-complex structure which are regular for simply covered disks with boundary in $L$. \eprop

It follows from \Cref{proposition:Jindep} that we can write $\tilde{n}(L, \a_0)= \tilde{n}(L, J, \a_0)$. 

Let us now consider some applications of \Cref{proposition:Jindep} to homological linking of Lagrangian submanifolds. We defer the proof of \Cref{proposition:Jindep} to the next section. 

\defi \label{definition:admissible} We say that a (non-monotone) Lagrangian torus $L \sub \R^4$ is \emph{admissible} if $\tilde{n}(L, \a_0) =1.$  \edefi

It is immediate that the class of admissible Lagrangian tori is closed under Hamiltonian isotopy. The next proposition shows that it contains ``most" product tori.

\prop \label{proposition:productexample} Consider the product torus $L(r,s)= \{(z_1, z_2) \in \C^2 \mid |z_1|=r, |z_2|=s \}.$ Assume without loss of generality that $0< r <s.$ Then $L(r, s)$ is admissible if $s/r \geq \sqrt{2}.$ The class $\a_0$ is represented by $[D^2 \tms *].$ \eprop

\pf Let us write $L= L(r,s)$. We first argue that $A_2(L) =\o(\a_0)= \pi r^2.$  By choosing each of the product factors as generators for $H_1(L; \Z),$ we get an identification $\Z \oplus \Z \simeq H_1(L; \Z) \simeq \pi_2(\R^4, L)$ sending $(1,0)$ onto $[D^2 \tms *]$ and $(0,1)$ onto $[* \tms D^2].$  Now, every Maslov $2$ class is of the form $(p, -p+1)$ for $p \in \Z,$ and so the areas of Maslov $2$ classes are of the form $\pi(p(r-s)+s).$  Using our assumption that $s \geq \sqrt{2} r,$ it is then easy to check that 
\eqs  A_2(L) = \pi r^2 = \op{min}_{p \in \Z} \lt( \{ \pi( p(r^2-s^2)+ s^2) \} \cap \R_{>0}  \rt). \eeqs  
It is now a standard fact that the standard complex structure $j$ is regular for all holomorphic disks with boundary in $L$, and that the boundary evaluation map has degree $1$; see \cite[Thm. 10.2]{cho} and \cite[Lem. 4]{auroux}. Hence $\tilde{n}(L, \a_0)=1$ and it follows that $L$ is admissible. \epf

\rmk \label{remark:workinprogress} It is a folklore conjecture that all non-monotone Lagrangian tori in $\R^4$ are Hamiltonian isotopic to product tori. In light of \Cref{proposition:productexample}, this would imply that the class of admissible tori contains ``most" examples of non-monotone Lagrangian tori in $\R^4$. \ermk

We are now in a position to prove the following result, which was stated as \Cref{theorem: nonmonotone} in the introduction. It extends \cite[Theorem 5.1]{evans-rizell} to the class of admissible tori in $\R^4.$ 

\thm \label{theorem: nonmonotone2} Let $L_1, L_2 \sub \R^4$ be disjoint Lagrangian tori. Suppose that $L_1$ is admissible. If $A_2(L_2) \geq A_2(L_1),$ then $[L_1]$ is the zero class in $H_n(\R^4-L_2; \Z).$ In other words, $L_1$ is homologically unlinked from $L_2$ (cf. \Cref{definition: ERlinking}). \ethm

\pf A routine modification of the proof of \Cref{proposition: sim1} (or equivalently \Cref{lemma:sim11} or \cite[Thm. 4.1]{evans-rizell}) shows that there exists a regular almost-complex structure $J^l$ with the property that the image of $\op{ev}(\a_0, J^l) \to \R^4$ misses $L_2.$ The boundary evaluation map has degree $1$ since $L_1$ is admissible. The rest of the argument is now identical to the proof of \cite[Thm. 5.1]{evans-rizell}. \epf

\subsection{Proof of \Cref{proposition:Jindep}} Let $J_0$ and $J_1$ be compatible almost-complex structures on $(\R^4, \o)$ which are regular for simply covered disks.  Let $\{J_t\}_{t \in [0,1]}$ be a generic homotopy of compatible almost-complex structures. Let $\mc{M}_{0,1}(\a, J_t) = \{ u: (D^2, \d D^2) \to (\R^4, L) \mid \ov{\d}_{J_t} u =0, u_*[D^2]= \a \} / \op{Aut}(D^2, 1).$

\lem \label{lemma: nonegmaslov} For $t \in [0,1]$, there does not exist a $J_t$-holomorphic disk $u: (D^2, \d D^2) \to (\R^4, L)$ such that $\mu(u)<0.$ \elem

\pf Assume for contradiction that such a disk exists. If $u$ is simply-covered, we get a contradiction due to the genericity of $\{J_t\}$ and the fact that $\op{ind}(u)= -1+ \mu(u) \leq -3.$ If $u$ is not simply-covered, there is a simply-covered $J_t$-holomorphic curve $v: (D^2, \d D^2) \to (\R^4, L)$ and a degree $d>1$ map $\phi: D^2 \to D^2$ such that $u= v \circ \phi.$ It follows that $\mu(v) = \mu(u)/d<0.$ Replacing $u$ with $v,$ we are back to considering the case where $u$ is simply-covered, which again gives a contradiction. \epf

\lem \label{lemma: nozeromaslov} For $t \in [0,1]$, suppose that $u: (D^2, \d D^2) \to (\R^4, L)$ is a $J_t$-holomorphic disk representing a class $\b \in \pi_2(\R^4, L).$ Let $\a_0$ be the $\mu$-infimal class of $L$. If $\mu(\b)=0,$ then $\o(\b) \geq \o(\a_0)= A_2(L).$ \elem

\pf By the long exact sequence of the pair $(\R^4, L),$ we have isomorphisms $\pi_2(\R^4, L) \simeq \pi_1(L) \simeq H_1(L; \Z)$. Hence we can view $\mu$ and $\o$ as elements of $\op{Hom}(H_1(L; \Z), \Z) \simeq H^1(L;\Z)$.

Let $\a_1 \in H_1(L; \Z)$ be the unique class such that $\mu(\a_1)=2, \o(\a_1) >\o(\a_0)$ and $\{\a_0, \a_1\}$ generates $H_1(L; \Z).$  To see that such a class exists and is unique, consider the preimage of $2$ under the group homomorphism $\mu: H_1(L; \Z) \to \Z.$ This is the intersection of a line in $H_1(L ; \R)$ with the lattice $H_1(L; \Z).$  It's not hard to check that any two adjacent lattice points on the line generate the entire lattice.  Now $\a_0$ has two adjacent lattice points, and $\a_1$ is the unique one satisfying the condition that $\o (\a_1)> \o (\a_0).$ 

Observe that $\o(\a_1) \geq 2\o(\a_0).$  Indeed, since the class $2\a_0- \a_1$ has Maslov index $2,$ it follows from the definition of $\a_0$ that either $0\geq \o(2\a_0-\a_1)$ or $\o(2 \a_0- \a_1)> \o(\a_0).$ The later inequality would contradict the fact that $\o(\a_1)> \o(\a_0).$ Hence $0 \geq \o(2\a_0-\a_1),$ which means that $\o(\a_1)\geq  2 \o(\a_0).$ 

Finally, it follows from the fact that $\{\a_0, \a_1\}$ generate $H_1(L;\Z)$ that all Maslov zero classes are of the form $n(\a_1-\a_0),$ for $n \in \Z.$  If $n \geq 1,$ then $\o(n(\a_1-\a_0)) = n \o (\a_1- \a_0) \geq n \o(\a_0) \geq A_2(L).$  If $n<0,$ then $\o(n(\a_1-\a_0))<0.$ Since classes of negative symplectic area do not support holomorphic disks, this proves the lemma. \epf

\prop \label{proposition: smooth} The moduli space $\mc{M}_{0,1}(\a_0, J_t)$ is a compact smooth manifold with boundary. Its boundary can be identified with $\mc{M}_{0,1}(\a_0, J_0) \sqcup \mc{M}_{0,1}(\a_0, J_1).$ \eprop

\pf It follows from the genericity of $\{J_t\}$ and the fact that $\a_0 \in \pi_2(\R^4, L)$ is primitive that $\mc{M}_1(\a_0, J_t)$ is a smooth manifold of finite dimension.  It remains to prove that it is compact. To this end, let $\{u_t\}$ be a sequence of $J_t$-holomorphic disks representing the class $\a_0$ which Gromov converge to a $J_{t_0}$-holomorphic stable holomorphic map $\bf{u}= (u^{\a})$ in the sense of \cite[Sec. 1.3]{frau-ziem}, for some $t_0 \in (0,1]$. Since $\pi_2(\R^4)=0,$ it follows that $u^{\a}: (D^2, \d D^2) \to (\R^4, L)$ for all $\a$ (i.e. there are no sphere bubbles). It follows from \Cref{lemma: nonegmaslov} that $\mu(u^{\a}) \geq 0$.  We claim that in fact $\mu(u^{\a}) \geq 2.$ If a Maslov $0$ disk occurred, it would follow by \Cref{lemma: nozeromaslov} that it would have maximal area and so there could be no other disks. But since we also have that $[u_t]= \sum_{\a} [u^{\a}]$ as classes in $\pi_2(\R^4, L),$ there must be other disks. Hence all $\mu(u^{\a})\geq 2.$ 

Since $\mu(u_t)=2,$ we conclude again from the fact that $[u_t]= \sum_{\a} [u^{\a}]$ that the stable map $\bf{u}=(u^{\a})$ consists of a single holomorphic disk of Maslov index $2$. This map must represent the class $\a_0$ since $\a_0=[u_t]= \sum_{\a} [u^{\a}]$. The proposition follows. \epf

\pf[Proof of \Cref{proposition:Jindep}] This follows from \Cref{proposition: smooth} and the fact that the degree is a cobordism invariant. \epf

\section{A construction of linked tori} \label{section: linkingconstruction}

The purpose of this section is to prove that the condition $A_2(L_2) \geq A_2(L_1)$ in the statement of \Cref{theorem: nonmonotone2} is sharp. More precisely, we prove the following proposition.

\prop \label{proposition: necessary1} Given real numbers $A^{(1)}> A^{(2)} >0$, there exists a pair of admissible, disjoint Lagrangian tori $L_1, L_2 \sub \R^4$ with $A_2(L_1)=A^{(1)} $ and $A_2(L_2)= A^{(2)}$, such that $[L_1] \in H_2(\R^4-L_2)$ is not the zero class.  \eprop

The proof of \Cref{proposition: necessary1} will proceed in three steps. We will first construct a pair of Lagrangian cylinders, such that one cylinders is ``threaded" through the other; cf. \Cref{lemma: linkedcylinder}. We will then ``close-up" these cylinders, thus obtaining a pair of Lagrangian tori. Finally, we will show that these tori satisfy the properties stated in \Cref{proposition: necessary1}. 

\subsection{Construction of linked cylinders}  \label{subsection:linkedcylinders} We begin with a definition. 

\defi \label{definition:lagrangiancylinder} A \emph{Lagrangian cylinder} in $\R^4$ is a Lagrangian submanifold which is diffeomorphic to $S^1 \tms \R$. A Lagrangian cylinder $L$ is said to be \emph{standard} if it is of the form \eq C(a,b;r)= \{ (x_1, y_2, x_2, y_2) \mid (x_1-a)^2+ (y_1-b)^2 = r^2,\, x_2=0 \},\eeq for some $(a, b, r) \in \R \tms \R \tms \R_{>0}$. We say finally that a Lagrangian cylinder is \emph{standard at infinity} if it agrees with a standard cylinder outside of some compact set. \edefi

Suppose that $A^{(1)}, A^{(2)}>0$ are positive real constants with $A^{(1)}> A^{(2)}$. Choose $r_1> r_2>0$ satisfying $\pi r_1^2= A^{(1)}$ and $\pi r_2^2= A^{(2)}.$ 

\lem \label{lemma: linkedcylinder} There exists a smooth embedding $\phi: \R/\Z \tms \R \to \R^4$ satisfying the following properties:
\begin{itemize}
\item[(i)] The image of $\phi$ is a Lagrangian cylinder.
\item[(ii)] We have $\phi(s, t)= (r_2 \cos 2 \pi s+ D, r_2 \sin 2 \pi s, 0, t)$ whenever $|t|>T$, for some fixed constants $T>0$ and $D\geq 2(r_1+r_2)$.
\item[(iii)] The curve $\phi(0,t)$ and the solid cylinder $\{(x_1, y_1, x_2, y_2) \mid x_1^2+ y_1^2 \leq r_1^2,\, x_2=0 \}$ with boundary $C(0,0; r_1)$ intersect transversally in a single point.
\item[(iv)] We have $\op{Im} \phi \cap C(0,0; r_1)= \emptyset$.  
\end{itemize}
\elem

\pf We prove \Cref{lemma: linkedcylinder} by describing a procedure to construct $\phi.$ We consider a map
\begin{align*}
\phi: \R/ \Z \tms \R &\to \R^4 \\
(s, t) &\mapsto (x(s,t), y(s,t), z(s,t), t).
\end{align*}

We wish to find sufficient conditions on the functions $x, y, z$ in order for $\phi$ to describe a parametrized Lagrangian embedding. Observe that the condition that $\phi$ be a Lagrangian embedding is equivalent to the equation
\eqs 0 = \o(\d_s \s, \d_t \s)= x_s y_t- x_t y_s+ z_s, \eeqs
or equivalently, 
\eq \label{equation: zs} z_s = -(x_s y_t - x_t y_s). \eeq

Set $\g_t(s)= (x(s,t), y(s,t)$. We can think of $\{\g_t\}_{t \in \R}$ as a $1$-parameter family of curves moving in $\R^2.$ Such a family is referred to as a ``Lagrangian movie" in \cite{sauvaget}. Observe from \eqref{equation: zs} that $z(s,t),$ and hence $\phi(s,t),$ is completely determined by $\g_t$ and $z(0,t).$ 

Since $z(0,t)= z(1,t),$ we must have
\begin{align} \label{equation:constraint}
0 = \int_0^1 z_s &= - \int_0^1 (x_s y_t-x_t y_s) = - \int_0^1 (x_s y_t + x_{ts} y ) = - \d_t \int_0^1 x_s y = \d_t \int_{\g_t} \l,
\end{align}
where we have used integration by parts in the third equality.

We obtain from the above computations the following necessary conditions for $\phi$ to determine a parametrized Lagrangian immersion: 
\begin{itemize}
\item[(i)] $\g_t$ is an immersion for all $t,$
\item[(ii)] $\d_t (\int_{\g_t} \l)=0.$
\end{itemize}

Observe that any family of immersions $\{\g_t\}$ which satisfies $\d_t (\int_{\g_t} \l)=0$ can be lifted to a Lagrangian immersion by specifying the map $t \mapsto z(0,t).$ Moreover, this map can be chosen arbitrarily. Observe finally that $\phi$ will be an embedding if, for all fixed $t \in \R$, the loop $\g_t$ has no self-intersections. (This condition is sufficient to ensure that $\phi$ is an embedding but is by no means necessary.)

It is now straightforward to construct $\phi$ satisfying the properties stated in \Cref{lemma: linkedcylinder}. One way of doing this is ensure that $\{ \g_t\}$ and $z(0,t)$ simultaneously satisfy the following conditions, where $D= 2(r_1+ r_2)$.  
\begin{itemize}
\item We have  $\g_t(s)= (x(s, t), y(s,t)) = (r_2 \cos 2 \pi s, r_2 \sin 2\pi s)$ and $z(s,t)= t$ for $|t| \leq 1$.  
\item There exists a constant $T \gg 0$ such that $\g_t(s) =(x(s, t), y(s,t)) = (r_2 \cos 2 \pi s + D, r_2 \sin 2\pi s)$ and $z(s,t)= 0$ for all $|t| \geq T$
\end{itemize}
For $1 \leq |t| \leq T,$ the movie $\{\g_t\}$ can be defined to simply translate the circle of radius $r_2$ centered at the origin to the circle of radius $r_2$ centered at the point $(D,0)$. For those values of $t$ such that $\g_t(s) \cap \{ x_1^2+y_1^2 = r_1^2\}$ is non-empty, one needs to choose $|z(0,t)|$ large enough so that $\phi$ does not intersect the solid cylinder $\{ x_1^2 + y_1^2 \leq r_1^2, x_2=0 \}$. 

The precise choice of $T$ is immaterial but can be taken to depend only on $D$. Observe that the condition $r_1>r_2$ is needed to ensure that $\op{Im} \phi \cap C(0,0, r_1)$ is empty. This completes the proof of \Cref{lemma: linkedcylinder}. \epf

\subsection{Closing-up the cylinders} \label{subsection:closeup} We fix $\phi, D, T$ as in \Cref{lemma: linkedcylinder}.

For $\delta_1, \delta_2 \gg T,$ we consider the truncations $C_1= C(0,0; r_1) \cap \{ |y_2| \leq \delta_1\}$ and $C_2= \op{Im} \phi \cap \{ |y_2| \leq \delta_2\}.$ We can assume that $\delta_2$ is large enough so that $C_2$ agrees with the standard cylinder $C(0,D; r_2)$ on the set $\{ \delta_2 -2 \leq |y_2| \leq \delta_2\}$.

The purpose of this section is to explain how to ``close-up" $C_1$ and $C_2$ by gluing to them suitable Lagrangian cylinders, in order to obtain Lagrangian tori $L_1$ and $L_2.$ These cylinders will be constructed in such a way that $L_1$ and $L_2$ are disjoint, both admissible, and satisfy $A_2(L_1)=\pi r_1^2=A^{(1)}$ and $A_2(L_2)= \pi r_2^2 = A^{(2)}$. 

We will only describe the construction of $L_2$ as the other case is similar and easier. 

Fix $\a \gg 1$ and $\delta \gg 1$ and consider an embedded curve $\g: [0, 5] \to \R^4$ with the following properties: 
\eqs \g(t) = \begin{cases}
 (D,0, 0,  \delta + t) &\text{for } t \in [0,1], \\
 (D,0, \a, (\delta+1)- (t-2)(2(\delta+1)) &\text{for } t \in [2,3], \\
 (D,0,0, -\delta+ (t-5) &\text{for }t \in [4,5].
 \end{cases} \eeqs
We also require that $ \g(t) \sub \{ y_2 \geq \delta+1 \}$ for $t \in [1,2]$ and that $\g(t) \sub  \{ y_2 \leq -(\delta+1)\}$ for $t \in [3,4].$

By the isotropic neighborhood theorem, we can construct a Lagrangian cylinder $C_{\g}$ in a neighborhood of $\g.$ We can assume that $C_{\g}$ has the property that
 $C_{\g} \cap \{ \delta \leq y_2 \leq \delta+1\} = \{ (x_1-D)^2+y_1^2= \e_2^2, x_2= 0, \delta \leq y_2 \leq \delta+1 \} \cup \{ (x_1-D)^2+y_1^2= \e_2^2, x_2= \a, \delta \leq y_2 \leq \delta+1 \}$, for some small constant $\e_2>0$. 

By rescaling and translating the cylinder $C_{\g}$ (thus possibly making $\a$ and $\delta$ larger), we can assume that $\e_2=r_2.$ 

If we set $\delta_2=\delta,$ then we can glue $C_{\g}$ to $C_2.$  We obtain a Lagrangian cylinder $L_2:= C_{\g} \cup C_2.$ 

The homology $H_1(L_2; \Z)$ is generated by the meridian $\s= \{(x_1-D)^2+ y_1^2= r_1^2, x_2=0, y_2= \delta\}$ and by a longitudinal curve $\tau$. We can assume that $\tau$ agrees with $\op{Im} \phi(0,-)$ on the set $\{|y_2| \leq \delta-1\}$.  

By choosing $\g$ appropriately, and choosing $\tau$ appropriately, we can ensure that the projection of $\tau$ to the $(x_1, y_1)$ plane has rotation number zero, and that the projection to the $(x_2,y_2)$ plane has rotation number $1.$ It follows that $\s$ and $\tau$ both have Maslov index $2.$ We can also ensure, by choosing $\a$ large enough, that the area of $\tau$ is arbitrarily large and in particular larger that $2 \pi r_2^2.$ Since all Maslov $2$ classes are of the form $\tau+ n(\tau-\s),$ for $n \in \Z,$ one readily verifies that $ A_2(L_2) = \pi r_2^2 = A^{(2)}.$ 

By the same argument, we can close up $C_1$ to obtain a Lagrangian torus $L_1$ with $A_2(L_1)= \pi r_1^2 = A^{(1)} >A^{(2)}$. Since we are free to choose $\delta_1 \gg \delta_2$, is evident from the construction that we can ensure $L_1$ and $L_2$ are disjoint.  

Since $\tau \sub L_2$ agrees with $\op{Im} \phi(0, -)$ for $\{|y_2| \leq \delta-1\}$, it follows from \Cref{lemma: linkedcylinder} (iii) and the construction of $L_1$ that there is a solid torus $S$ with $\d S= L_1$ such that $\tau$ intersects $S$ transversally in a single point. It follows by an argument analogous to the proof of \Cref{lemma: pi1link} that the map $\pi_1(L_2) \to \pi_1(\R^4-L_1)$ has nontrivial image. We conclude in light of \Cref{lemma: pi1link} that $[L_1] \in H_2(\R^4-L_2; \Z)$ is not the zero class, i.e. $L_1$ is homologically linked with $L_2$.

\subsection{Admissibility} It remains to show that $L_1$ and $L_2$ are admissible. We will again only prove this for $L_2$ as the argument for $L_1$ is essentially the same. 

Let $A= \{ \delta \leq y_2 \leq \delta+1\} \sub \R^4.$ Observe that $L_2 \cap A=  \{ (x_1-D)^2+y_1^2 = r_2^2, x_2=0, \delta \leq y_2 \leq \delta + 1 \} \cup \{ (x_1-D)^2+y_1^2= r_2^2, x_2=\a, \delta \leq y_2 \leq \delta + 1 \}.$ 

Let $j$ be the standard integrable complex structure on $\R^4.$ Let $\tilde{j}$ be a small perturbation having the following properties.
\begin{itemize}
\item[(i)] $\tilde{j}$ is standard at infinity. 
\item[(ii)] $\tilde{j}$ agrees with $j$ on $A$. 
\item[(iii)] Any simply-covered $\tilde{j}$-holomorphic disk with boundary in $L_2$ having a point mapped into $\R^4-A$ is regular.
\end{itemize}

Observe that there are two families of embedded $\tilde{j}$-holomorphic disks parametrized by $s \in [0,1]$ which are of the form
\eqs \phi_s^1(x, y) = (x+D, y, 0, \delta+s) \eeqs
and 
\eqs \phi_s^2(x, y)= (x+D, y, \a, \delta+s). \eeqs

\lem Every simply-covered $\tilde{j}$-holomorphic disk either has a point mapped into $\R^4-A$ or belongs to one of the families $\{ \phi_s^i\},$ up to reparametrization. \elem

\pf Let $u: (D^2, \d D^2) \to (\R^4, L_2)$ be $\tilde{j}$-holomorphic and assume that $\op{Im}(u) \sub A.$ We write $u=(u_1, u_2)$ where $u_i$ is the projection onto the $(x_i, y_i)$ plane for $i=1,2.$ Since $\tilde{j}$ is standard on $A,$ it follows that the $u_i$ are ordinary holomorphic functions.

Observe that $u_2( \d D^2) \sub \{ x_2=0, \delta \leq y_2 \leq \delta+1\} \cup \{ x_2=\a, \delta \leq y_2 \leq \delta+1 \}.$ Since $u(\d D^2)$ is connected, it must be entirely contained in either of these two intervals. Let us assume that $u(\d D^2) \sub  \{ x_2=0, \delta \leq y_2 \leq \delta+1\}$ as the other case can be treated in the same way. 

We claim that in fact $\op{Im}(u_2) \sub \{ x_2=0, \delta \leq y_2 \leq \delta+1 \}.$ Assume for contradiction that this is not the case. Writing $u_2=(u_2^{x_2}, u_2^{y_2}),$ there exists $(x_0, y_0) \in \op{Im}(u_2)$ with the property that $|x_0|>0$ and $|u_2^{x_2}(x,y) |\leq |x_0|$ for all $(x,y) \in D^2.$ Hence there exists a point $p \in \op{Int}(D^2)$ with $u_2(p)= (x_0, y_0).$ This contradicts the open mapping theorem. \epf

\cor The almost-complex structure $\tilde{j}$ is regular for all simply-covered holomorphic disks with boundary in $L_2.$ \ecor
\pf It follows by automatic transversality that the $\tilde{j}$-holomorphic disks $\phi_s^i$ are all regular since they are embedded. We also know that $\tilde{j}$ is regular for all simply-covered holomorphic disks which have a point mapped in $\R^4-A.$ It follows from the lemma that there are no other simply-covered holomorphic disks. \epf

Consider the class $\a_0= [ \phi_s^i] \in \pi_2(\R^4, L_2)$. It follows from the construction of \Cref{subsection:closeup} that $\a_0$ is the $\mu$-infimal class of $L_2$; cf. \Cref{definition:muinfimal}. Let $\op{ev}(\a_0, J): \mc{M}_{0,1}(\a_0, \tilde{j}) \to L_2 \sub \R^4$ be the boundary evaluation map. 

\lem The degree of $\op{ev}(\a_0, J)$ is $1.$ \elem

\pf Choose $p \in \{ (x_1-D)^2+ y_1^2= r_2^2, x_2=0, \delta< y_2 < \delta+1 \} \sub L_2 \cap A.$ Suppose $p \in u(\d D^2).$ We claim that $u$ is a reparametrization of a curve in one of the families $\{ \phi_s^i\}.$ 

As before, write $u=(u_1, u_2)$ and observe that $u$ is an ordinary holomorphic function in $u^{-1}(A).$ Suppose first that $\op{Im}(u_2) \cap A  \sub \{ x_2=0, \delta \leq y_2 \leq \delta+1\} \cup \{ x_2=\a, \delta \leq y_2 \leq \delta+1\}.$ Since $u_2$ is an ordinary holomorphic function in $u^{-1}(A)$, it follows that $u_2$ is a constant function. From this, we can easily conclude that $u$ is a reparametrization of a curve in one of the family $\{ \phi_s^1\}.$ 

Suppose now that $\op{Im}(u_2) \cap A$ contains a point $(x_0, y_0)$ such that $x_0 \notin \{ 0, \a\}.$ It follows by the open mapping theorem that in fact $\{ 0 \leq x_2 \leq \a, \delta \leq y_2 \leq \delta+1\} \sub \op{Im}(u_2).$ 

But observe now that $\a= \int_{u^{-1}(A)} u_2^* \o_2 \leq \int_{u^{-1}(A)} u^* \o \leq \int_{D^2} u^* \o.$ This contradicts the fact that $\o(u)= \pi r_2^2$ since $\a \gg 1$ (and in particular, we were free to assume when choosing $\a$ that $\a > \pi r_2^2$).  \epf

We have shown that $\op{ev}(\a_0, \tilde{j})$ is a degree 1 map, where $\a_0$ is the $\mu$-infimal class and $\tilde{j}$ is regular for simply-covered disks with boundary in $L_2$. We conclude that $L_2$ is admissible; cf. \Cref{definition:admissible}. This concludes the proof of \Cref{proposition: necessary1}.

\section{Quantitative unlinking} \label{section:quantitativeunlink}

\subsection{Motivation} Consider a Lagrangian torus $L_1 \sub \R^4$. Let $\g \sub \R^4-L_1$ be a non-contractible embedded loop and let $\mc{U} \sub \R^4-L_1$ be a tubular neighborhood of $\g$. Now suppose that $L_2 \sub \mc{U}$ is another Lagrangian torus. Observe that if the map $\pi_1(L_2) \to \pi_1(\mc{U})= \Z$ has nontrivial image, then it follows from \Cref{lemma: pi1link} that $L_1$ is homologically linked with $L_2$. This means that any obstruction to linking $L_1$ with $L_2$ automatically gives an obstruction to embedding $L_2$ into $\mc{U}$ in a homologically essential way. 

In light of the results of \Cref{section: enumerativeinvariant}, where we showed that the linking behavior of tori is sensitive to the enumerative invariants $n(L, \a)$, one expects that such invariants could also be used to obstruct embeddings of tori into certain subdomains of $\R^4$.  We will see one instance of this in \Cref{proposition: cliffordnonsqueezing}. 

The discussion of this section can be fit neatly into the framework of symplectic capacities. Thus most of the results we present will be deduced from the existence of a certain symplectic capacity, which is a slight variant on a construction of Cieliebak and Mohnke in \cite[p.\ 2]{cmaudin}. 

The existence of this capacity can be deduced from a theorem of Charette \cite[Thm. 3.1]{charette} which was proved by Floer theoretic methods; cf. \Cref{remark:charette}. We will however present a self-contained proof which is closer in spirit to the arguments of the previous sections. 

This section is logically independent from the rest of the paper and the results may already be known to experts. Nevertheless, we feel that it serves a useful purpose in highlighting some connections between the study of linking and certain classical questions in symplectic topology.

\subsection{A symplectic capacity} As usual, we identify $\C^2$ with $\R^4$ by letting $(z_1, z_2)= (x_1+ iy_1, x_2+iy_2).$  Let us consider the polydisk 
\eqs \mc{P}(a, b) = \{(z_1, z_2) \in \C^2 \mid |z_1| < a, |z_2|< b \}. \eeqs

We can view $\mc{P}(a,b)$ both as an open symplectic manifold and as a symplectic subdomain of $(\R^4, \o)$. 

\prop \label{proposition: polydiskembedding} For $0<a \leq b$, let $L \sub \mc{P}(a,b) \sub \R^4$ be a Lagrangian torus. Then $A_2(L) \leq \pi a^2 .$ \eprop

\pf Suppose for contradiction that $A_2(L)> \pi a^2.$ As in \Cref{subsection:neckstretch}, the polydisk $\mc{P}(a,b)$ naturally embeds as a symplectic subdomain of $(S^2 \tms S^2, \o_1 \oplus \o_2)$, where $\int_{S^2} \o_1= \pi a^2$ and $\int_{S^2} \o_2 = \pi b^2.$  We write $D_{\infty}= S^2 \tms S^2 - \mc{P}(a,b)=  S^2 \tms \{ \infty \} \cup \{ \infty \} \tms S^2.$

Let $J$ be an almost-complex structure on $S^2 \tms S^2$ compatible with $\o_1 \oplus \o_2.$ For some $p \in S^2,$ let $ \a= [S^2 \tms p] \in H_2( S^2 \tms S^2 ; \Z).$  We consider the moduli space \eqs \mc{M}_1(\a, J):= \{ u: S^2 \to S^2 \tms S^2 \mid \ov{\d}_J u=0, u_*[D^2]= \a \} / \op{Aut}(S^2, x),\eeqs for some $x \in S^2$, and the evaluation map
\begin{align*}
\op{ev}(\a, J): \mc{M}_1(\a, J) &\to S^2 \tms S^2 \\
[u] &\mapsto u(x). 
\end{align*}
 
As in \Cref{subsection:neckstretch} and in \cite{igr}, we let $N(L) \sub S^2 \tms S^2$ be a Weinstein embedding and let $\mc{U}$ be an open neighborhood of $D_{\infty}$ such that $N(L) \cap \mc{U}= \emptyset$. Let $J^0$ be a compatible almost-complex structure on $S^2 \tms S^2$ which agrees with the standard integrable complex structure on $\mc{U}$. We now construct a sequence of compatible almost-complex structures $\{J^k\}_{k=1}^{\infty}$ by stretching the neck along $S^*_{1, g}\bb{T}^2 \sub N(L)$, where $g$ is a suitably rescaled flat metric on $\bb{T}^2$.

A well-known theorem of Gromov implies that $\op{ev}(\a, J^k)$ is a degree $1$ map (and in fact the $J^k$-holomorphic spheres in the class $\a$ form a foliation of $S^2 \tms S^2$). It follows from the SFT compactness theorem that there is an infinite sequence $u_k$ of $J^k$-holomorphic spheres which converge to a building $\mathbf{u}.$ Let $J^{\infty}$ be the almost-complex structure on $S^2 \tms S^2- L$ which results from the neck-stretching procedure. By choosing $J^0$ appropriately, or equivalently by simultaneously perturbing the $J^k$ in the complement of $N(L) \cup \mc{U}$, we can assume that $J^{\infty}$ is regular for simply covered punctured holomorphic curves. 

It follows by elementary topological considerations that the building $\mathbf{u}$ must have at least two $J^{\infty}$-holomorphic planes in $S^2 \tms S^2-L$. A routine modification of the proofs of \Cref{lemma: next} and \Cref{lemma:posindex} shows that the components of $\mathbf{u}$ satisfy the following two properties.
\begin{itemize}
\item[(i)] The sum of the Fredholm indices of all the components of $\mathbf{u}$ which map into $S^2 \tms S^2 -L$ is at most $2$. 
\item[(ii)] Every component $u^{\tau}$ of $\mathbf{u}$ has non-negative Fredholm index. If moreover $u^{\tau}$ is a plane, then $\op{ind}(u^{\tau}) \geq 1$ with equality if and only if $u^{\tau}$ is simply-covered and the compactification $\ov{u}^{\tau}$ has Maslov index $2$. 
\end{itemize}

We conclude that there are exactly two simply-covered planes of Fredholm index $1$. It follows by positivity of intersection (since $J^{\infty}$ is standard in $\mc{U}$) that only one such plane can intersect $D_{\infty}.$ Let $v$ be the plane which does not intersect $D_{\infty}.$ We can think of $v$ as a plane inside $\mc{P}(a,b) \sub \R^4.$ It follows by \Cref{lemma:maslovchern} that $\mu(\ov{v})=2.$ But this implies that $\o(\ov{v}) \geq A_2(L),$ contradicting our assumption that $\o(\ov{v}) \leq \o(u_k) = \pi a^2 < A_2(L).$ \epf

\rmk \label{remark:charette} As noted above, \Cref{proposition: polydiskembedding} can also be deduced from work of Charette \cite{charette, charette1}. In fact, one can prove the stronger statement that $A_2(L) \leq d(L)$, where $d(L)$ is the displacement energy of $L$. This follows by combining \cite[Thm. 3.1]{charette} and the fact that there are no holomorphic disks of Maslov index strictly less than 2 for a generic almost-complex structure. \ermk

Following Cieliebak and Mohnke \cite[p. 2]{cmaudin}, \Cref{proposition: polydiskembedding} can interpreted in terms of a symplectic capacity. 

\defi For any domain $U \sub \R^4,$ we define a symplectic capacity $c_{L,2}$ as follows:
\eqs c_{L,2} (U):= \op{sup} \{ A_2(L) \mid L \sub U \; \text{embedded Lagrangian torus} \} \in [0, \infty]. \eeqs

\edefi

It is clear that this capacity is well-defined and nonzero on any non-empty domain, since we can always embed a Clifford torus with sufficiently small monotonicity factor.  

\prop The capacity $c_{L,2}$ satisfies the following properties:

\begin{tabular}{lp{10.5cm}}
(Monotonicity)&We have $c_{L,2}(U') \leq c_{L,2}(U)$ if $\mc{U}' \sub \mc{U}$. \\ 
(Conformality)&Given any real constant $r>0,$ we have $c_{L,2}(rU)= |r|^2 c_{L,2}(U).$\\
(Invariance)&If $\phi$ is a Hamiltonian isotopy, then $c_{L,2}(U)= c_{L,2}(\phi(U))$.\\
(Nontriviality)&We have $0<c_{L,2}(B^4(1))$ and $c_{L,2}(Z^4(1))< \infty,$ where we write $Z^4(1)= \R^2 \tms B^4(1) \sub \R^4.$
\end{tabular}
\eprop

\pf The first two properties are immediate from the definition. The invariance property follows immediately from the fact that $A_2(L)$ is invariant under Hamiltonian isotopy. The fact that $0<c_L^2(B^4(1))$ is also clear since we may embed a Clifford torus with monotonicity factor $\pi^2/8$ inside the unit ball.  Finally, the fact that $c_L^2(Z^4(1))< \infty$ is an immediate consequence of \Cref{proposition: polydiskembedding}. \epf

In fact, \Cref{proposition: polydiskembedding} implies that $c_{L,2}(Z^4(1)) \leq \pi.$  It's clear that $\pi \leq c_{L,2}(Z^4(1))$ since we may embed a Clifford torus of factor $\frac{\pi}{2}(1-\e)^2$ for any $\e>0.$ It follows that $c_{L,2}(Z^4(1))=\pi.$   

Similarly, let $0 < a \leq b$ be as in \Cref{proposition: polydiskembedding} and observe that we can embed a Clifford torus of factor $\frac{\pi}{2}(a-\e)^2$ inside $\mc{P}(a,b).$ It follows that $\pi a^2 \leq c_{L,2}(\mc{P}(a,b)).$ The proposition now implies that $c_{L,2}(\mc{P}(a,b)) \leq \pi a^2$ which means that $c_{L,2}(\mc{P}(a,b)) = \pi a^2.$ 

We are thus led to the following result: 

\prop \label{proposition: easyembedding} Consider the product torus $L(r,s) = \{ (z_1, z_2) \mid |z_1|=r, |z_2|=s \}$ for $0< r \leq s.$  If $s \geq \sqrt{2} r,$ and $r>a,$ then $L(r,s)$ cannot be embedded by a Hamiltonian isotopy into the polydisk $\mc{P}(a,b)$.  It also follows that the polydisk $\mc{P}(r,s)$ cannot be embedded by a Hamiltonian isotopy into the polydisk $\mc{P}(a,b)$. \eprop

\pf[Proof of \Cref{proposition: easyembedding}] By the same argument as in the proof of \Cref{proposition:productexample}, one can show that the condition $s \geq \sqrt{2} r$ implies that $A_2(L) = \pi r^2.$ Since we observed above that $c_{L,2}(\mc{P}(a,b))= \pi a^2$, it follows from our assumption that $r>a$ and from the definition of $c_{L,2}$ that $L(r,s)$ cannot be embedded in $\mc{P}(a,b)$. 

The fact that $\mc{P}(r,s)$ cannot be embedded into the polydisk $\mc{P}(a,b)$ follows from the monotonicity and invariance properties of the capacity. \epf

\rmk \label{remark:cs} \Cref{proposition: easyembedding} can also be deduced from work of Chekanov and Schlenk, who proved in \cite[Sec.\ 2.1]{chek-schlenk} that the displacement energy of $L(r,s)$ for $r \leq s$ is $\pi r^2$. \ermk

\ex Suppose that $r=1$ and $s=3/2.$ It follows from the proposition that the torus $L(r,s)$ cannot be embedded by a Hamiltonian isotopy into $\mc{P}(1-\e, 1-\e)$ for any $\e>0.$ Observe that there exists a class in $H_1(L(r,s); \Z)$ (or Maslov index $-1$) having symplectic area $-2\pi+ \frac{9}{4} \pi = \frac{\pi}{4}.$ Hence the capacity $c_{L}$ defined in \cite{cmaudin} does not a priori rule out the existence of such an embedding.  Of course, the embedding can easily be ruled out by the work of Chekanov and Schlenk mentioned in \Cref{remark:cs}. \eex

\subsection{Quantitative non-linking} In \Cref{proposition: easyembedding}, we gave obstructions to Hamiltonian embeddings of Lagrangian tori into certain polydisks in terms of the invariant $A_2(L).$ The purpose of this section is to establish an obstruction to Hamiltonian embeddings of Lagrangian tori into certain subdomains of $\R^4$ which depends on the enumerative invariants $n(L, \a_i)$ considered in \Cref{section: enumerativeinvariant}.

Consider the domain $\mc{D}_{\e} = \{ 1-\e < |z_1|< 1+\e \} \tms \{ |z_2|< \e \} \sub \C^2$ for any $\e \in (0,1/2).$ 

\prop \label{proposition: topologicalnonsqueezing} Let $L \hookrightarrow \mc{D}_{\e}$ be an embedded Lagrangian torus. Suppose that there exist classes $\a_1, \a_2 \in \pi_2(\R^4,L)$ satisfying properties (i), (ii), (iii) of \Cref{proposition: essential}. Then the natural map $\pi_1(L) \to \pi_1(\mc{D}_{\e}) = \Z$ is trivial.  \eprop

\cor \label{proposition: cliffordnonsqueezing} Let $L_{Cl} \hookrightarrow \mc{D}_{\e}$ be an embedded Clifford torus. Then the induced map $\pi_1(L_{Cl}) \to \pi_1( \mc{D}_{\e})= \Z$ is trivial.  Note that we do not make any assumptions about the monotonicity factor of $L_{Cl}.$ \ecor

\pf It follows from \Cref{proposition: polydiskembedding} that $A_2(L) \leq \pi \e^2.$  Now assume that $\pi_1( L ) \to \pi_1(\mc{D}_{\e})$ is nontrivial.  Without loss of generality, we can assume that $\a_1$ has nontrivial image in $\pi_1(\mc{D}_{\e}).$ This implies that $\a_1$ has nontrivial image in $\pi_1(\mc{D}_{\e}) \simeq \pi_1(\C^* \tms \C) \simeq \pi_1(\C^*)$.

Let $J$ be a compatible almost-complex structure which is regular for simply-covered disks with boundary in $L$ and is obtained by perturbing the standard integrable complex structure $j$ in the interior of $\mc{D}_{\e}$.  Since $n(L, \a_1) =1$, it follows that $\mc{M}(\a_1, J)$ is non-empty.  Let $u \in \mc{M}(\a_1, J).$  Let $u_1:= \pi_1 \circ u,$ where $\pi_1: \C^2 \to \C$ is the projection onto the first factor.

Observe that $u_1$ is holomorphic on $u_1^{-1}( \{ |z_1|< 1-\e\})$ since $J$ is standard on this domain.  Since $[\d u]= \op{Im}(\a_1)  \in \pi_1( \C^*)$ is nontrivial, it follows that $\d u_1$ has nontrivial winding number.  Hence $0 \in \op{Im} u_1.$  It now follows by the open mapping theorem that $\op{Im} u_1 \cap \{ |z_1|< 1-\e\} =  \{ |z_1|< 1-\e\}.$  But 
\eq \pi (1-\e)^2 = \int_{u_1^{-1}(\{|z_1|< 1-\e)\} } u_1^* \o_1 \leq \int_{u_1^{-1}(\{|z_1|< 1-\e)\}} u^* \o \leq \int_{D^2} u^* \o. \eeq  The middle inequality uses the fact that $u$ is holomorphic on $ u_1^{-1}(\{ |z_1| < 1-\e \}).$ This contradicts the fact that $A_2(L) \leq \pi \e^2.$  \epf

In contrast to \Cref{proposition: cliffordnonsqueezing}, there is no obstruction to squeezing Chekanov tori.

\prop There exists an embedded Chekanov torus $L_{Ch} \to \mc{D}_{\e}$ such that the induced map $\pi_1(L_{Ch}) \to \pi_1(\mc{D}_{\e})$ is surjective. \eprop
\pf Choose a simple closed curve $\g \sub \mc{D}_{\e}$ which represents a nontrivial class in $\pi_1(\mc{D}_{\e}).$  The desired claim now follows from \Cref{corollary: chekanovsqueeze}. \epf

\section{Closing remarks} We end this paper by briefly discussing to what extent our methods are limited to dimension 4. We also highlight some possible directions for further research.

\subsection{The role of dimension 4} As a general rule, all results in this paper which rely on the analysis of holomorphic planes are expected to fail in dimensions greater than $4$. This applies in particular to the results on smooth unlinking in \Cref{section:monotonelink} and \Cref{section: enumerativeinvariant}.  These results make essential use of the intersection theory of \cite{siefring} and \cite{siefringwendl} and of index positivity properties (cf. \Cref{lemma:posindex}), neither of which are available in dimension greater than 4. The importance of the intersection theory is partly hidden from view in our paper since it enters into the proof of \Cref{proposition: igrciting}, which we obtained as a consequence of arguments in \cite{igr}.

On the other hand, the methods of \Cref{section: non-monotone} do work in higher dimensions and should in principle allow one to prove homological linking results for non-monotone tori in all dimensions.  However, these results would get progressively weaker as the dimension increases, in the sense that the corresponding class of ``admissible tori" would get smaller. This is essentially because one needs to prevent the appearance of disks of Maslov index $[2-n, 0]$ in order to prove an analog of \Cref{proposition: smooth} for tori in $\R^{2n}$. We also remark that the class of admissible tori in $\R^4$ is plausibly very large; cf. \Cref{remark:workinprogress}. This is unlikely to be true in higher dimensions. 

Regarding the constructions of \Cref{section: linkingconstruction}, it is certainly possible to construct Lagrangians in all dimensions as lifts of lower dimensional projections (see, for instance, the technique of ``Lagrangian suspension" \cite[3.1E]{polt-book}). However, the 4-dimensional case is particularly easy to visualize and to work with, because the projections are just closed curves in $\R^4$ and the area constraint \eqref{equation:constraint} takes a very simple form. This allows us to effectively ``see" what types of movies are possible.

\subsection{Smooth unlinking of admissible Lagrangian tori} \label{subsection:planesvsdisks} One could hope to improve the results of \Cref{section: non-monotone} from homological unlinking to smooth unlinking. A natural approach would be to work with holomorphic planes rather than disks.  If one tries to implement this approach, one runs into the difficulty that holomorphic planes can degenerate into a priori complicated buildings which could potentially have certain non-regular and multiply-covered components. In contrast, holomorphic disks in an exact symplectic manifold can only degenerate into disks. Although the analysis of planes appears more complicated, there is reason to hope that it could be tractable. In particular, one could hope to take advantage of the intersection theory of \cite{siefring} and \cite{siefringwendl}, and of the many useful results contained in \cite[Sec. 3 \& 4]{igr}. 

\subsection{Lagrangian unlinking} In light of the classification result of Dimitroglou Rizell, Ivrii and Goodman \cite{igr} (see \Cref{theorem: igr}), one might hope to upgrade the results of \Cref{section:monotonelink} and \Cref{section: enumerativeinvariant} from smooth unlinking to Lagrangian unlinking. For example, one might hope to show that two Clifford tori can always be pulled apart from each other through Lagrangian tori; cf. \Cref{theorem:clifforddisjoint}.  

If one follows the strategy of \cite[Sec. 6]{igr}, the key step in constructing Lagrangian isotopies is to extend the embedded solid tori constructed using holomorphic planes to a symplectic embedding of $S^1 \tms D^2 \tms (-\e, \e)$. To achieve this, one needs to ensure that the symplectic disks which foliate the solid torus have trivial monodromy. For a single torus, this can be achieved using the so-called ``inflation procedure", at the cost of modifying the Lagrangian; cf. \cite[Sec. 6]{igr}. However, if there are two or more Lagrangians, a naive application of the inflation procedure might cause them to intersect. 

\subsection{Linking in high dimensions} It could be interesting to study the connection between the enumerative invariants of the type considered in \Cref{section: enumerativeinvariant} and linking of tori in higher dimensions.  What can one say about linking of Clifford tori in high dimensions? Auroux \cite{auroux} has constructed infinite families of monotone Lagrangian tori in $\R^{2n}$ for $n \geq 3$ which are not of Clifford or Chekanov type, and which are distinguished by enumerative invariants analogous to those considered in \Cref{section: enumerativeinvariant}. What sort of linking behavior do these tori display?  Note that in high dimensions, one only expects to obtain results about homological linking by directly analyzing moduli spaces of holomorphic curves. However, it might still be possible in favorable circumstances to promote such results to statements about smooth linking, using techniques of high dimensional differential topology.  

\subsection{Local linking} In \cite[Thm. 1.1B]{eliash}, Eliashberg proved a ``local unknottedness result" which states that any Lagrangian cylinder in $\R^4$ which is standard at infinity is Hamiltonian isotopic to a standard cylinder (see \Cref{definition:lagrangiancylinder}). In this spirit, one could also try to prove ``local unlinking" results. One expects that a cylinder of radius $r$ should be smoothly unlinked from any cylinder of radius $R \geq r$. Given a configuration of $N$ disjoint cylinders in $\R^4$ which are standard at infinity, one also expects this configuration to be smoothly isotopic to some standard model which depends only on how the various components are homologically linked. Finally, the monodromy issues mentioned above do not occur for cylinders, so it should be possible to upgrade a smooth isotopy to a Lagrangian isotopy using techniques from \cite[Sec. 6]{igr}. 

As in \Cref{section:monotonelink}, a first step in proving such statements would be to show that the cylinders under consideration occur as the boundary of an embedded solid cylinder which is foliated by holomorphic planes.  One way to do this would be to start with a family of planes near infinity (where the cylinder is standard) and to argue that this family is open and closed and hence extends to the whole cylinder.  This can be done using the theory of \cite{siefring} and \cite{siefringwendl}, although the analysis is not completely straightforward due to the non-compactness of the domain. (Alternatively, one could also compactify the situation and turn the cylinder into a torus. This would solve the compactness issue but one would lose the monotonicity of the cylinder).

\section{Appendix}

The purpose of this appendix is to briefly collect some definitions and computations from the theory of punctured pseudoholomorphic curves which are needed in \Cref{section:monotonelink} of this paper. We will assume that the reader is familiar with the basics of this theory, as outlined for instance in \cite{wendlautomatic}. The definitions and notation below are intended to be consistent with \cite{wendlautomatic}.

\subsection{The index formula for punctured pseudoholomorphic curves} \label{section: basicdefis} 

Let $(W, J)$ be an almost-complex manifold with cylindrical ends of the form $(-\infty, 0] \tms M_-$ and $[0, \infty) \tms M_+$ (the manifolds $M_+$ and $M_-$ are are allowed to be disconnected or empty).

Given a punctured Riemann surface $\dot{\S}= (\S- \G^+- \G^-),$ let $\bf{c}$ denote an assignment of a family of Reeb orbits to each puncture $z \in \G^+ \cup \G^-.$ The Fredholm index of a $J$-holomorphic curve $u: \dot{\S} \to (W,J)$ with asymptotic orbits determined by $\bf{c}$ is shown \cite[eq. (1.1)]{wendlautomatic} to satisfy the formula 
\eq \label{equation: indexformula} \op{ind}(u ; \bf{c})= (n-3) \chi(\dot{\S}) + 2c_1^{\Phi}(u^*TW) + \mu^{\Phi}(u; \bf{c}).\eeq

Let us briefly recall how the terms appearing in \eqref{equation: indexformula} are defined. Here $c_1^{\Phi}(u^*TW)$ is the relative first Chern number of $u^*TW$ with respect to a trivialization $\Phi$ near the punctures of $\dot{\S}.$  It counts the number of zeros of a generic section of $u^*TW \wedge u^*TW$ which is constant and nonzero near the punctures with respect to the trivialization induced by $\Phi.$  

Our definition of $\mu^{\Phi}$ follows \cite[Sec. 3.2]{wendlautomatic}.  Given a $T$-periodic orbit $c_{\g}$ of the contact manifold $M_{\pm}$, there is an associated asymptotic operator 
\begin{align*}
\bf{A}_{\g}: \G(x^*\xi) &\to \G(x^*\xi)\\
\bf{A}_{\g}&= -J( \grad_t v - T \grad_v R_{\a}), 
\end{align*}
where $x: S^1 \to M_{\pm}$ is a parametrization of $c_{\g}$ satisfying $\dot{x} = c_{\g} = TR_{\g}.$ 

The Conley-Zehnder index of a non-degenerate asymptotic operator is defined as in \cite[Definition 3.30]{wendlsft}.  If $\bf{A}_{\g}$ is a degenerate asymptotic operator, then the operator $(\bf{A}+ \delta \op{Id})$ is non-degenerate provided that $\delta \notin \s(\bf{A}_{\g}).$  One can show that $\mu_{CZ}^{\Phi}(\g_z \pm \delta):= \mu_{CZ}^{\Phi}(\bf{A}_{\g_z} \pm \delta \op{Id})$ is well-defined provided that $\delta$ is sufficiently small. 

We will always be considering punctured holomorphic curves $\S$ with \textit{unconstrained} ends $\bf{c}$ which are allowed to move in a Morse-Bott manifold, and possibly with boundary component $\d \S$ contained in some Lagrangian $L \sub \R^{2n}.$  In this case, we have for $\delta>0$
\eq \mu^{\Phi}(u; \bf{c})= \sum_{z \in \G^+} \mu_{CZ}^{\Phi}(\g_z-\delta)- \sum_{z \in \G^-} \mu_{CZ}^{\Phi}(\g_z+ \delta) + \mu( \d \S). \eeq

\subsection{Some index computations} Let us now specialize to the setting of \Cref{section:monotonelink}. We will follow throughout this section the notation introduced in \Cref{subsection:neckstretch}. In particular, recall that $S^*_{1, g_i}\bb{T}^2$ is a contact manifold with coordinates $(\t_1, \t_2, \ov{\t})$ and contact form $\a_i$ for $i=1,2$. We have trivializations $\Phi^i= \{ \d_t, R_{\a}, X= \sin \ov{\t} \d_{\t_1} - \cos \ov{\t} \d_{\t_2}, \d_{\ov{\t}} \}$ of the tangent bundle of $\R \tms S^*_{1, g_i}\bb{T}^2$.  With respect to $\Phi^i$, the almost complex structure $J_{cyl}$ takes the form 
\eq J_{cyl}= \begin{pmatrix}
0 & -1& 0 &0 \\
1 & 0 &0 &0 \\
0 & 0 & 0 &-1 \\
0 &0 &1 & 0 
\end{pmatrix}. \eeq

We wish to compute the Conley-Zehnder index of Reeb orbits in $S^*_{1, g_i}\bb{T}^2$ with respect to trivializations $\{X, \d_{\ov{\t}} \}$ of $\xi= \ker \a_i.$  

\lem \label{lemma: CZindex}  Let $c_{\g}$ be a Reeb orbit of period $T$ with $c_{\g} = R_{\g} = \cos \ov{\t} \d_{\t_1} + \sin \ov{\t} \d_{\t_2}.$   Then $\mu_{CZ}^{\Phi}(\g - \delta)= 1$ while $\mu_{CZ}^{\Phi}(\g+ \delta) = 0.$ 
\elem

\pf We have: 
\begin{align*}
\bf{A}_{\g} + \delta \op{Id} &= \begin{pmatrix}
 0 & 1 \\
-1 & 0 \end{pmatrix}
\lt(
\begin{pmatrix} 
\d_t & 0 \\
0 & \d_t 
\end{pmatrix} 
-T 
\begin{pmatrix}
0 & -1\\
0 & 0 
\end{pmatrix} 
\rt)
+ \begin{pmatrix}
\delta & 0 \\
0 & \delta
\end{pmatrix}
=
\begin{pmatrix}
\delta & \d_t \\
- \d_t & -T + \delta 
\end{pmatrix}.
\end{align*}

The constant vectors $X, \d_{\ov{\t}}$ are eigenvectors with eigenvalues $\delta$ and $-T+\delta$ respectively, each of winding number number zero.  It follows from \cite[Theorem 3.36]{wendlsft} that $\mu_{CZ}^{\Phi}(\g+\delta) =0$ while $\mu_{CZ}^{\Phi}(\g-\delta)= 1.$ \epf

In the context of \Cref{section:monotonelink}, we are considering punctured holomorphic curves mapping into almost-complex manifolds with cylindrical ends diffeomorphic to $(-\infty, 0] \tms S^*_{1, g_i} \bb{T}^2$ or $[0, \infty) \tms S^*_{1, g_i} \bb{T}^2$ and endowed with the almost-complex structure $J_{cyl}$. 

In this situation, it follows from the index formula \eqref{equation: indexformula} and from \Cref{lemma: CZindex} that the index of a pseudoholomorphic curve $u$ satisfies the formula:
\begin{align} \label{equation: explicitindex}
\op{ind}(u) &= -\chi(u) + \# \{ \text{positive punctures of } u \} -0 + 2c_1^{\Phi}(u)  \\
&=-2+ \# \{\text{all punctures of } u \} + \# \{ \text{positive punctures of } u \} + 2c_1^{\Phi}(u). \nonumber
\end{align}

\subsection{The Maslov index and relative Chern number} The purpose of this section is to sketch a proof of \Cref{lemma:maslovchern}, which we restate here for the reader's convenience. We will follow throughout this argument the notation of \Cref{subsection:neckstretch}.

\begin{lemma}[cf. \Cref{lemma:maslovchern}] \label{lemma:maslovchernappendix} For $i=1,2$, let $u_i: \C \to (S^2 \tms S^2 -L_i)$ be a $J$-holomorphic plane where $J= J_{\infty}^l$ or $J=J_k^{\infty}$. Let $v: \C \to S^2 \tms S^2 -L_1-L_2$ be a $J$-holomorphic plane for $J= J_{\infty}^{\infty}$. Then $2c_1^{\Phi}(u_i)= \mu(\ov{u_i})$ and $2c_1^{\Phi}(v) = \mu(\ov{v})$. \end{lemma}

\pf Throughout this argument, let $u$ stand for $u_1, u_2,v.$  

There are many related but distinct notions in the literature which go by the name of ``Maslov index". In \cite[Chap. 2]{mcduff-sal-intro}, one considers a Maslov index for loops of symplectic m
atrices and a Maslov index for loops of Lagrangian subspaces in $(\R^{2n}, \o)$. Let us denote the former by $m_s$ and the latter by $m_l$. These two indices are related as follows. If $\g: S^1 \to \mc{L}(n)$ is a loop of Lagrangian subspaces and $\s: S^1 \to \op{Sp}(2n)$ is a loop of symplectic matrices such that $\s(t) \g(0)= \g(t)$, then one has (see \cite[Thm. 2.3.7]{mcduff-sal-intro}) 
\eq \label{equation:maslovrelation} 2m_s(\s)= m_l(\g). \eeq 

The Maslov class $\mu(\ov{u})$ is defined as follows. Observe that there is a homotopically unique trivialization $\tau$ of $u^*T(S^2 \tms S^2)$. There is also a path $\g: \d D^2 \to u^*T(S^2 \tms S^2)|_{\d D^2}$ of Lagrangian subspaces determined by $L_i$, which can be viewed as a path $\g: \d D^2 \to \mc{L}(2)$ with respect to the trivialization $\tau$. We then have \eq \label{equation:maslovrelation2} \mu(\ov{u})= m_l(\g). \eeq

Now observe that $\Phi= \{ \d_t, R_{\a_i}, X, \d_{\ov{\t}} \}$ extends to a $C^0$ trivialization of $\ov{u}^*T(S^2 \tms S^2)|_{\d D^2}$ with the property that the subframe $\{ R_{\a_i}, X\}$ is tangent to $L_i$. Let $\Phi(t)= \Phi|_{u(\g(t))}$ and let $\s$ be a loop of linear maps such that $\Phi(t)= \s(t) \Phi(0)$. Then $\s$ can be viewed as a loop of symplectic matrices with respect to the homotopically unique trivialization of $\ov{u}^*T(S^2 \tms S^2)$. It follows from \eqref{equation:maslovrelation} and \eqref{equation:maslovrelation2} that $2m_s(\s)= \mu(\ov{u})$. 

It only remains to relate $m_s(\s)$ to $c_1^{\Phi}(u)$, i.e. we wish to relate the winding number of $\Phi$ with respect to the homotopically unique trivialization of $\ov{u}^*T(S^2 \tms S^2)$ with the count of zeros of a generic section which is constant near the punctures with respect to $\Phi$. It can be shown by standard arguments (see \cite[Sec. 2.7]{mcduff-sal-intro}) that these counts are equal. \epf

\subsection{Notions of energy}  \label{subsection:energyappendix} As we noted in \Cref{remark:sft}, the proof of \Cref{proposition: sim1} uses a version of SFT compactness for ``stretching the neck" in a manifold with a negative cylindrical end. To the author's knowledge, a proof of this precise version of the theorem has not appeared in the literature. However, closely related statements are proved in \cite{sftcompactness} and the arguments there go through in our setting with routine modifications. 

In order to apply the approach of \cite{sftcompactness} for proving SFT compactness, one needs to control certain energies of punctured holomorphic curves. In \Cref{proposition: sim1}, one controls the symplectic area of a sequence of $J_{\infty}^l$-holomorphic planes $u_l$. The purpose of this section is to argue that our control on the symplectic areas of these planes automatically gives us control on their energies in the sense needed for applying the arguments of \cite{sftcompactness}. The following arguments are mainly drawn from \cite{sftcompactness} but are included for completeness. 

Let us briefly recall the setting of \Cref{proposition: sim1}. We are considering $J_{\infty}^l$-holomorphic planes $u_l$ for $l>0$ in the symplectic manifold $(S^2 \tms S^2 -L_1, \o)$. 

Let us write $V:= S^*_{1,g_1}\bb{T}^2 \sub S^2 \tms S^2-L_1$, where $S^*_{1,g_1}\bb{T}^2$ is a unit circle bundle as defined in \Cref{subsection:neckstretch}. Recall that the symplectic manifold $S^2 \tms S^2 -L_1$ has a negative end of the form $$((-\infty, 1] \tms V, d(e^t \a) ),$$ where we let $\a$ be the contact form on $V$ which was defined in \Cref{subsection:neckstretch} and denoted there by $\a_1$.  It will be convenient to write $$S^2 \tms S^2-L_1= E_- \cup_V M_+,$$ where we let $E_-= ((-\infty, 0] \tms V)$ and $M_+= S^2 \tms S^2 - ((-\infty, 0) \tms V)$.

Let $v: \dot{\S} \to S^2 \tms S^2-L_1$ be a punctured $J_{\infty}^l$-holomorphic curve.  Let $\pi: (-\infty, 1] \tms V \to (-\infty, 1]$ be the projection. A standard argument using the maximum principle implies that $\pi \circ v$ has no critical points. Hence $v^{-1}(t_0 \tms V) \sub \S$ is a manifold with boundary for all $t_0 \in (-\infty, 1]$. Let us write $v= (s, g)$ on $v^{-1}((-\infty, 1] \tms V) . $

We first define the notions of energy considered in \cite{sftcompactness} which are relevant in our setting. We begin with the the so-called $\a$-energy.  To this end, let $\mc{C}$ be the set of all functions $\phi: \R_- \to \R_+$ such that $\int \phi =1$.  

\defi[see (23) in \cite{sftcompactness}] We define the \emph{$\a$-energy} of a $J_{\infty}^l$-holomorphic curve $v$ as follows: $$E_{\a}(v):= \sup_{\phi \in \mc{C}} \int_{v^{-1}(E_-)} (\phi \circ s) ds \wedge g^* \a.$$ \edefi

Next, we consider the $\o$-energy: 

\defi[see (22) in \cite{sftcompactness}] We define the \emph{$\o$-energy} of a $J_{\infty}^l$-holomorphic curve $v$ as follows: $$E_{\o}(v):= \int_{v^{-1}(E_-)} v^* d\a + \int_{v^{-1}(M_+)} v^* \o.$$ \edefi

\rmk The term ``$\o$-energy" is potentially confusing since it does not coincide with the symplectic area, but we have retained it to be consistent with \cite{sftcompactness}. \ermk

Our goal is now to prove \Cref{corollary:maincor}. This says that the symplectic area of the planes $u_l$ considered in \Cref{proposition: sim1} controls -- up to a constant factor which is independent of $l$ -- the $\a$- and $\o$-energies.  For the remainder of this section, we will abuse notation by dropping the subscript $l$ from our notation and writing $u= (t, f)$ on $u^{-1}((-\infty, 1] \tms V)$. The reader may verify that the following inequalities are independent of $l$.

Let's begin by analyzing the $\o$ energy. 

\lem \label{lemma:dalpha} We have $\int_{u^{-1}(E_-)} du^*(\a) \leq  \int_{u^{-1}(0 \tms V)} u^*(\a)$. \elem
\pf  Let $\g$ be the Reeb orbit to which $u$ is asymptotic.  Let $\mc{A}(\g):= \int_{\g} \a$ be the action of $\g$. By Stokes' theorem, we have: 
\eq \int_{u^{-1}(E_-)} du^*(\a) = \int_{u^{-1}(0 \tms V)} u^*(\a)- \mc{A}(\g)  \leq \int_{u^{-1}(0 \tms V)} u^*(\a). \eeq 
 \epf

Hence, it is enough to control $\int_{u^{-1}(0 \tms V)} u^*(\a)$.   This is the content of the next lemma

\lem[cf.\ Lem.\ 9.2 of \cite{sftcompactness}] \label{lemma:keyslicebound} We have $\int_{u^{-1}(0 \tms V)} u^*(\a) \leq \frac{1+e}{e-1} \omega(u)$.  \elem

\pf Let $$C_1= \int_{u^{-1}([0,1] \tms V)} u^{*} d(e^t \a)= \int_{u^{-1}(1 \tms V)} e u^{*} \a -  \int_{u^{-1}(0 \tms V)} u^{*} \a.$$  

Also let $$C_2= \int_{u^{-1}([0,1] \tms V)} u^{*} d\a= \int_{u^{-1}(1 \tms V)} u^{*} \a -  \int_{u^{-1}(0 \tms V)} u^{*}\a.$$

Hence $C_1-eC_2= (e-1)  \int_{u^{-1}(0 \tms V)} u^{*} \a$.  Note that we evidently have $C_1 \leq \omega(\ov{u})$. 

We claim that $D \leq \omega({u})$. Indeed, since $u$ is holomorphic, it follows that $$C_2= \int u^* d\a \leq \int u^*( e^t(d\a)) \leq \int u^*(e^t( dt \wedge \a + d\a)) \leq \omega({u}).$$

Hence \eq \int_{u^{-1}(0 \tms V)} u^{*} \a = \frac{C_1-eC_2}{e-1} \leq \frac{1+e}{e-1} \omega(u). \eeq \epf

By combining the above two lemmas and appealing to the definition of $E_{\o}$, we obtain the following corollary:

\cor We have $E_{\o}(u) \leq (1+  \frac{1+e}{e-1}) \omega(u)$. \ecor

Let us now analyze the $\a$-energy. Recall that $\mc{C}$ is the set of functions $\phi: \R_- \to \R_+$ such that $\int \phi=1$. 

Given $\phi \in \mc{C}$, we let $\psi(s)= \int_{-\infty}^s \phi(t) dt$. 

\lem \label{lemma:alphaenergy} We have $E_{\a}(u) \leq  \int_{u^{-1}(0 \tms V)} f^* \a$. \elem

\pf We write: 
\begin{align*} \text{sup}_{\phi \in \mc{C}} \int_{u^{-1}(E_-)} (\phi \circ t) dt \wedge f^* \a &=  \int_{u^{-1}(E_-)} (d (\psi(t) f^* \a)- \psi(t) f^* d\a) \\
&= \int_{u^{-1}(0 \tms V)} f^* \a - \int_{u^{-1}(E_-)} \psi(t) f^* d\a.
\end{align*}

The second term is always non-positive since $u$ is holomorphic, so this proves the claim.\epf 

\cor We have $E_{\a}(u) \leq  \frac{1+e}{e-1} \omega(u)$. \ecor

\pf This follows by combining \Cref{lemma:keyslicebound} and \Cref{lemma:alphaenergy} . \epf

\cor \label{corollary:maincor} We have $E_{\a}(u)+ E_{\o}(u) \leq (1+ 2\frac{1+e}{e-1}) \omega(u)$. \ecor

\end{document}